\documentclass[a4paper,11pt]{article}
\usepackage{a4wide}

\usepackage[T1]{fontenc}
\usepackage[utf8]{inputenc}
\usepackage{amsmath}
\usepackage{mathtools}
\usepackage{amsthm}
\usepackage{amsfonts}
\usepackage{amssymb}
\usepackage{empheq}
\usepackage{graphicx}
\usepackage{subcaption}
\usepackage{xcolor}
\usepackage{tikz}
\usepackage{listofitems}
\usepackage{multirow}

\usetikzlibrary{decorations.pathmorphing,patterns}
\tikzstyle{mynode}=[thick,draw=blue,fill=blue!20,circle,minimum size=12,inner sep=0pt]

\usepackage{hyperref}

\graphicspath{{Figures/}}

\begin{document}

\title{Comparison of neural and spline representations\\ for physics-informed learning}
\author{R\'egis Duvigneau\footnote{Inria, 2004 route des lucioles, 06902 Sophia-Antipolis, France - regis.duvigneau@inria.fr}}
\date{Universit\'e C\^ote d'Azur, Inria, CNRS, LJAD}

\maketitle

\begin{abstract}
Physics-informed machine learning as recently emerged as a new paradigm for the resolution of partial differential equations, including functionalities for parametric modeling and inverse design, using neural networks to represent the solution fields. However, several difficulties have been reported concerning the training of such networks, resulting in slow convergence, low accurate solutions and sensitivity to the choice of hyper-parameters. In order to better understand the mechanisms underlying these difficulties, the present work proposes a rigourous comparison of neural networks and B-Splines representations within the same physics-informed learning framework. Four problems with closed-form solutions are numerically studied to establish a rigorous analysis of the accuracy and efficiency of the two representations, in terms of asymptotic convergence rates, sampling requirements, loss function minimization and arithmetic precision needed.
\end{abstract}

\section{Introduction}

The current context of energy transition makes the use of multi-physics simulation and design methodologies critical across all industrial sectors (energy, transportation, building, chemicals, etc.). Continuous advances in numerical simulation methods over the last thirty years have enabled the development of design workflows incorporating numerous cutting-edge methods, such as geometric modeling, mesh generation, numerical resolution of partial differential equations (PDEs), optimization. Although these techniques are now well established, their wide adoption in industry has yet to be achieved, primarily for two reasons: the \emph{complexity} of implementing these workflows and the \emph{computational cost} of running them. These observations motivate the development of new methods that are simpler to implement and less costly to run.

Recently, the emergence of machine learning techniques has enable the construction of approximate models based on data, whether from experimental measurements or pre-existing simulation results, which have been found particularly effective for many applications, when large databases are available.
Many architectures have been investigated, from historical multilayer perceptrons~\cite{anderson1995introduction} to convolutional networks~\cite{cnn_cfd,Lee_You_2019}, graph neural networks~\cite{GNS_20} and, more recently, transformers~\cite{geneva2022transformers,wu2402transolver} including attention mechanisms and neural operators~\cite{li2020fourier, lu2021learning} mapping functions.

Among these, the present work focuses on techniques that make direct use of physical models based on PDEs during the training phase, by incorporating the PDE residuals in the loss function, which allows to address data scarcity but also to solve inverse problems. These approaches are known, amongst other names, as Physics-Informed Neural Networks (PINNs)~\cite{raissiPhysicsinformedNeuralNetworks2019}, deep Ritz method~\cite{yu2018deep} or deep Galerkin method~\cite{sirignano2018dgm}.
They have enjoyed considerable success in recent years, leading to a remarkable number of publications, numerous extensions and applications across all disciplines, as demonstrated in the review article~\cite{cuomoScientificMachineLearning2022}.
However, these physics-informed learning approaches cannot be considered as general-purpose solvers yet, such as finite-element method. Indeed, several studies have pointed out strong limitations, in particular the sensitivity with respect to the choice of hyper-parameters such as network architecture, optimizer, loss weights, etc.~\cite{krishnapriyan2021characterizing,wangExpertGuideTraining2023,xu2025fp64needrethinkingfailure}.

The objective of this work is to study more in depth the role of the neural representation itself in these difficulties, by comparing the use of neural networks and B-Splines in a similar physics-informed machine learning framework. B-Splines are a standard representation in geometric modeling~\cite{DeBoor_78,Farin_89} and are also used for PDE solving in isogeometric analysis methods~\cite{Cottrell_etal_09}. Recently, they have also been employed in machine learning in conjunction with neural networks~\cite{saidaoui2024deep,sun2023pisl,wang2025physics}, mainly to facilitate geometry handling. In a very recent work~\cite{canali2026trainable}, B-Splines have been used directly to parameterize the solution in a physics-informed machine learning formulation. The present work is based on similar methodological ingredients, but is focused on the analysis of the approximation properties of the two representations.

The general framework of physics-informed learning is presented in Section~\ref{sec:piml}, irrespective of the representation chosen for the solution. The specific features of neural networks and B-Splines are described in Section~\ref{sec:nn} and~\ref{sec:spline}. Their properties are examined in Section~\ref{sec:prop}. A set of four problems with closed-form solutions is considered in Section~\ref{sec:result}, allowing a rigorous comparison of errors obtained by numerical resolution. Finally, the results are discussed in Section~\ref{sec:disc}, before conclusion.

\section{Physics-informed machine learning}
\label{sec:piml}

\subsection{Principle of residuals minimization}

Physics-informed machine learning aims to simulate physical phenomena by minimizing an error estimate based on model residuals, typically formalized as PDEs. Although the approach is clearly inspired by artificial intelligence paradigms, the use of PDE residuals makes it similar to standard resolution methods like finite-element method.

For the sake of simplicity, we consider a generic problem to present the methodology: the physical rules are expressed as a one-dimensional PDE, first-order in time and second-order in space, associated with Dirichlet and Neumann boundary conditions and initial condition. The extension to more complex governing equations will be straightforward. Thus, the problem writes:

\begin{align}
\begin{cases}
    \partial_t{u} + \mathcal{N}[u,\partial_x{u},\partial_{xx}{u}](x,t) = 0 & (x,t) \in \Omega \times ]0,T] \\
    u(x^D,t) = g^D(x^D,t) & (x^D,t) \in \partial \Omega^D \times ]0,T]\\
    \partial_n{u}(x^N,t) = g^N(x^N,t) & (x^N,t) \in \partial \Omega^N \times ]0,T]\\
    u(x,0) = g^I(x) & x \in \Omega ,
\end{cases}
\label{eq:pb}
\end{align}

where $u$ is the unknown field, $\Omega$ an open domain, $(x,t)$ the space-time coordinates and $\mathcal{N}$ denotes a differential operator. $g^D$, $g^N$ and $g^I$ define respectively the Dirichlet, Neumann and initial conditions.

The resolution of the problem is based on a parametric representation of the solution $u(\theta; x,t)$, where $\theta$ is a set of parameters to be calibrated,  also denoted as \emph{degrees of freedom}. The precise nature of the representation will be examined in forthcoming sections, most of the methodological ingredients of this approach being agnostic w.r.t. the chosen representation. The calibration of the parameters $\theta$, denoted as \emph{training}, simply consists in minimizing an estimate of the error $\mathcal{L}(\theta)$, denoted as \emph{loss function}, which aggregates errors regarding PDE residuals and initial / boundary conditions:
\begin{equation}
    \mathcal{L}(\theta) = \mathcal{L}_{PDE}(\theta) + \mathcal{L}_{I}(\theta) + \mathcal{L}_{D}(\theta) + \mathcal{L}_{N}(\theta) .
\end{equation}
$\mathcal{L}_{PDE}$ is the loss term related to the PDE residual:
\begin{equation}
    \label{eq:pdeloss}
    \mathcal{L}_{PDE}(\theta) = \frac{1}{N_{PDE}} \sum_{k=1}^{N_{PDE}} \left \{ \mathcal{R}(\theta;x_k,t_k) \right \}^2 ,
\end{equation}
with:
\begin{equation}
    \mathcal{R}(\theta;x_k,t_k) = \partial_t{u}(\theta;x_k,t_k) + \mathcal{N}[u,\partial_x{u},\partial_{xx}{u}](\theta;x_k,t_k) .
\end{equation}
The points $(x_k,t_k)$, denoted as \emph{collocation points}, correspond to a space-time sampling of the computational domain, where the PDE has to be fulfilled. Thus, $\mathcal{L}_{PDE}$ approximates the square of the $L_2$-norm of the residuals.
Similarly, other loss terms concern initial and boundary conditions:
\begin{equation}
    \mathcal{L}_{I}(\theta) = \frac{1}{N_{I}} \sum_{k=1}^{N_{I}} \left \{ u(\theta;x_k,0) - g^I(x_k,0) \right \}^2  ,
\end{equation}
\begin{equation}
    \mathcal{L}_{D}(\theta) = \frac{1}{N_{D}} \sum_{k=1}^{N_{D}} \left \{ u(\theta;x_k^D,t_k) - g^D(x_k^D,t_k) \right \}^2  ,
\end{equation}
\begin{equation}
    \mathcal{L}_{N}(\theta) = \frac{1}{N_{N}} \sum_{k=1}^{N_{N}} \left \{ \partial_n{u}(\theta;x_k^N,t_k) - g^N(x_k^N,t_k) \right \}^2 ,
\end{equation}
with dedicated samplings of initial and boundary subdomains.

As can be observed, this approach is agnostic with regard to the representation of the solution. The only requirement is the ability to evaluate the space and time derivatives present in the loss terms.

\subsection{Training procedure}

The minimization of the loss function is the core of the approach, which determines the practical implementation and the accuracy and robustness of the method. Unfortunately, this may be a difficult optimization problem, characterized by a high-dimensional set of parameters $\theta$ and, as pointed out by several studies~\cite{krishnapriyan2021characterizing,wangExpertGuideTraining2023}, a possibly non-convex and highly curved loss function lanscape. Due to the large dimension, the task is usually achieved using descent methods, updating the parameters according to:
\begin{equation}
\theta^{(p+1)} = \theta^{(p)} - \rho^{(p)} \mathcal{G}(\nabla_\theta \mathcal{L}(\theta^{(p)})) ,
\end{equation}
where $\rho$ is the step length, also denoted as \emph{learning rate}, and $\nabla_\theta \mathcal{L}$ the gradient of the loss function. $\mathcal{G}$ is an operator that defines the algorithm. In most studies, the optimization is achieved in two distinct phases: a first one based on a stochastic first-order descent algorithm allowing a robust start, like ADAM~\cite{kingma2017adammethodstochasticoptimization}, followed by a second-order quasi-Newton algorithm yielding a faster convergence, like L-BFGS~\cite{liu1989limited}. Nevertheless, it remains a tedious optimization problem and several approaches have been proposed in the literature to improve the convergence~\cite{cuomoScientificMachineLearning2022}.

The evaluation of the gradient of the loss function necessitates the differentiation of the different loss terms:
\begin{equation}
    \nabla_\theta \mathcal{L}(\theta) = \nabla_\theta \mathcal{L}_{PDE}(\theta) + \nabla_\theta \mathcal{L}_{I}(\theta) + \nabla_\theta \mathcal{L}_{D}(\theta) + \nabla_\theta \mathcal{L}_{N}(\theta) ,
\end{equation}
with:
\begin{equation}
    \nabla_\theta \mathcal{L}_{PDE}(\theta) = \frac{2}{N_{PDE}} \sum_{k=1}^{N_{PDE}} \left \{ \mathcal{R}(\theta;x_k,t_k) \right \} \, \partial_\theta \mathcal{R}(\theta;x_k,t_k) ,
\end{equation}
\begin{equation}
    \nabla_\theta \mathcal{L}_{I}(\theta) = \frac{2}{N_{I}} \sum_{k=1}^{N_{I}} \left \{ u(x_k,0) - g^I(x_k) \right \} \, \partial_\theta u(\theta;x_k,t_k) ,
\end{equation}
\begin{equation}
    \nabla_\theta \mathcal{L}_{D}(\theta) = \frac{2}{N_{D}} \sum_{k=1}^{N_{D}} \left \{ u(\theta;x_k^D,t_k) - g^D(x_k^D,t_k) \right \} \, \partial_\theta u(\theta;x_k,t_k) ,
\end{equation}
\begin{equation}
    \nabla_\theta \mathcal{L}_{N}(\theta) = \frac{2}{N_{N}} \sum_{k=1}^{N_{N}} \left \{ \partial_n{u}(\theta;x_k^N,t_k) - g^N(x_k^N,t_k) \right \} \, \partial_\theta u(\theta;x_k,t_k)  .
\end{equation}
As can be noticed, it is required to evaluate the derivative of the solution with respect to the parameters $\partial_\theta u$, as well as the derivative of the residuals $\partial_\theta \mathcal{R}$. As the latter involves space and time derivatives of the solution, this results in a higher-order cross derivative evaluation.

In the next sections, two possible representations for the solution are presented, focusing on their respective characteristics.

\section{Neural representation}
\label{sec:nn}

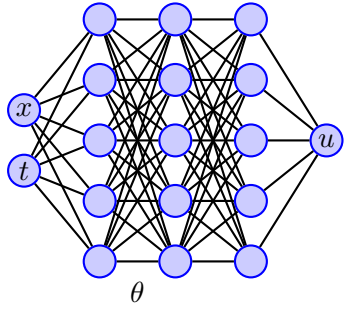
\begin{figure}[!ht]
\centering
\begin{tikzpicture}[x=1cm,y=0.8cm]

\node at (2.5,-2.5)  {$\theta$};
  \readlist\Nnod{2,5,5,5,1} 
  \foreachitem \N \in \Nnod{ 
    \foreach \i [evaluate={\x=\Ncnt; \y=\N/2-\i+0.5; \prev=int(\Ncnt-1);}] in {1,...,\N}{ 
      \node[mynode] (N\Ncnt-\i) at (\x,\y) {};
      \ifnum\Ncnt>1 
        \foreach \j in {1,...,\Nnod[\prev]}{ 
          \draw[-,thick] (N\prev-\j) -- (N\Ncnt-\i); 
        }
      \fi 
    }
  }
\node at (1,0.5)  {$x$};
\node at (1,-0.5)  {$t$};
\node at (5,0)  {$u$};

\end{tikzpicture}
\caption{Multi-layer perceptron with 3 hidden layers of 5 neurons.}
\label{fig:mlp}
\end{figure}

The Multi-Layer Perceptron (MLP) network is used in the PINNs methodology~\cite{raissiPhysicsinformedNeuralNetworks2019} to represent the solution of the problem. It is composed of a set of L layers: one input layer, $L-2$ intermediate hidden layers and one output layer. Each intermediate layer $l$ contains a set of $N_l$ neurons, each of them being connected to all neurons of the previous layer. The number of neurons in the input layer corresponds to the number of space and time coordinates, whereas the number of neurons in the output layer is the number of unknown fields, as depicted in Fig.~(\ref{fig:mlp}).

Thus, the values associated with the neurons of the first layer are simply the space-time coordinates $(x,t)$ of a given collocation point. The value of the neuron in the last layer predicts the solution $u(\theta;x,t)$. Each neuron of the internadiate layers is a mathematical operator that applies a non-linear \emph{activation function} $\alpha$ to the weighted sum of its own inputs plus a bias factor. Thus, the value of the $i$th neuron of the layer $l$ writes:
\begin{equation}
f_i^l = \alpha(\sum_{j=1}^{N_{l-1}} w^l_{ij} f_j^{l-1} + b_i^l)
\end{equation}
where $N_{l-1}$ is the number of neurons in the $l-1$th layer, $w^l_{ij}$ the weight of the connection between the neuron $i$ and the neuron $j$ of the layer $l-1$ and $b_i^l$ is the bias. Hyperbolic tangent is usually employed as activation function for PINNs, because of its smoothness and non-linearity. No activation is usually applied to the last layer, then the predicted solution is just a linear combination of the $N_{L-1}$ values obtained in the last but one layer. The use of the hidden layers can therefore be interpreted as the construction of $N_{L-1}$ suitable functions to represent the solution.
As summary, the output $u(\theta;x,t)$ is obtained by propagating forward the input values $(x,t)$ through the network layers, resulting in a composition of the activation functions.
Such a network is characterized by its weights and biases, which constitute the set of parameters $\theta$ to be calibrated during the training phase.

The evaluation of the space and time derivatives consists in differentiating the output of the network with respect to its inputs. According to the forward evaluation achieved by composition of the activation function, the most efficient approach to estimate the derivatives is based on a backpropagation algorithm~\cite{Rojas1996}, usually implemented using Automatic Differentiation (AD) techniques~\cite{JMLR:v18:17-468}. The derivatives of the loss function with respect to the parameters $\theta$ (network weights and biases), required for the descent algorithm, are evaluated in the same way.

\section{Spline representation}
\label{sec:spline}

B-Splines are considered as a standard in Computer-Aided Design~\cite{DeBoor_78,Farin_89}, since they allow to represent exactly a broad class of curves / surfaces / volumes commonly encountered in engineering and permit a very intuitive definition and modification of geometrical objects through the handling of \emph{control points}. In the present context, B-Spline functions are not in charge of describing a geometry, but representing the solution of a PDE problem.

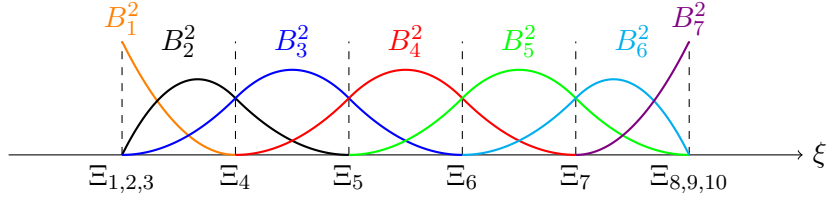
\begin{figure}[!ht]
    \begin{center}
    \begin{tikzpicture}[scale = 1.5]
    \draw[->] (-1,0) -- (6,0) node[right] {$\xi$};
    \draw[dashed] (0,0) -- (0,1) ;
    \node[below] at (0,0) {$\Xi_{1,2,3}$} ;
    \draw[dashed]  (1,0) -- (1,1) ;
    \node[below] at (1,0) {$\Xi_{4}$} ;
    \draw[dashed]  (2,0) -- (2,1) ;
    \node[below] at (2,0) {$\Xi_{5}$} ;
    \draw[dashed]  (3,0) -- (3,1) ;
    \node[below] at (3,0) {$\Xi_{6}$} ;
    \draw[dashed]  (4,0) -- (4,1) ;
    \node[below] at (4,0) {$\Xi_{7}$} ;
    \draw[dashed]  (5,0) -- (5,1) ;
    \node[below] at (5,0) {$\Xi_{8,9,10}$} ;

    \draw[domain=0:1,smooth,variable=\x,orange,thick] plot ({\x},{(1-\x)*(1-\x)});
    \node[orange] at (0,1.2) {$B^2_1$} ;

    \draw[domain=0:1,smooth,variable=\x,black,thick] plot ({\x},{2*\x-3/2*\x*\x});
    \draw[domain=1:2,smooth,variable=\x,black,thick] plot ({\x},{0.5*(2-\x)*(2-\x)});
    \node[black] at (0.5,1) {$B^2_2$} ;

    \draw[domain=0:1,smooth,variable=\x,blue,thick] plot ({\x},{0.5*\x*\x});
    \draw[domain=1:2,smooth,variable=\x,blue,thick] plot ({\x},{-0.5+\x-(\x-1)*(\x-1)});
    \draw[domain=2:3,smooth,variable=\x,blue,thick] plot ({\x},{5/2-\x+0.5*(\x-2)*(\x-2)});
    \node[blue] at (1.5,1) {$B^2_3$} ;

    \draw[domain=1:2,smooth,variable=\x,red,thick] plot ({\x},{0.5*(\x-1)*(\x-1)});
    \draw[domain=2:3,smooth,variable=\x,red,thick] plot ({\x},{-0.5+(\x-1)-(\x-2)*(\x-2)});
    \draw[domain=3:4,smooth,variable=\x,red,thick] plot ({\x},{5/2-(\x-1)+0.5*(\x-3)*(\x-3)});
    \node[red] at (2.5,1) {$B^2_4$} ;

    \draw[domain=2:3,smooth,variable=\x,green,thick] plot ({(\x},{0.5*(\x-2)*(\x-2)});
    \draw[domain=3:4,smooth,variable=\x,green,thick] plot ({\x},{-0.5+(\x-2)-(\x-3)*(\x-3)});
    \draw[domain=4:5,smooth,variable=\x,green,thick] plot ({\x},{5/2-(\x-2)+0.5*(\x-4)*(\x-4)});
    \node[green] at (3.5,1) {$B^2_5$} ;

    \draw[domain=3:4,smooth,variable=\x,cyan,thick] plot ({\x},{0.5*(\x-3)*(\x-3)});
    \draw[domain=4:5,smooth,variable=\x,cyan,thick] plot ({\x},{-16+10*(\x-1)-3/2*(\x-1)*(\x-1)});
    \node[cyan] at (4.5,1) {$B^2_6$} ;

    \draw[domain=4:5,smooth,variable=\x,violet,thick] plot ({\x},{(\x-4)*(\x-4)});
    \node[violet] at (5,1.2) {$B^2_7$} ;

    \end{tikzpicture}
    \end{center}
\caption{Seven quadratic B-Spline basis functions.}
\label{fig:bs_basis}
\end{figure}

The B-Spline basis functions are defined using a so-called \emph{knot} vector  $(\Xi_1, \dots, \Xi_l)\in \mathbb{R}^{l}$, which consists of  $l$ nondecreasing real numbers. This knot vector defines a discretization of the \emph{parametric domain} $[\Xi_1 , \Xi_l]$.
B-Spline basis functions $(B^p_{i})_{i=1, \cdots, n}$ are polynomials and defined recursively in the parametric domain as~\cite{DeBoor_78}:
\begin{equation}
  B^0_{i}(\xi) = \begin{cases} 1 &  \text{ if } \Xi_i \leq \xi  < \Xi_{i+1}\\
    0 & \text{ otherwise}
  \end{cases}
  \label{eq:nurbs0}
\end{equation}
\begin{equation}
  B^p_{i}(\xi) = \frac{\xi- \Xi_i}{\Xi_{i+p} - \Xi_i} B^{p-1}_{i}(\xi) +
  \frac{\Xi_{i+p+1}- \xi}{\Xi_{i+p+1}-\Xi_{i+1}} B^{p-1}_{i+1}(\xi).
  \label{eq:nurbs_eval}
\end{equation}
Note that the quotient 0/0 is assumed to be zero.
The degree of the functions $p$, the number of knots $l$ and functions $n$ are related by $l=n+p+1$~\cite{DeBoor_78}.
Clamped knot vectors, i.e. knot vectors with first and last knots of multiplicity $p+1$,
are usually used for representations of degree $p$ to impose interpolation and tangency conditions at both extremities~\cite{DeBoor_78}. This property can be exploited to impose initial, Dirichlet and Neumann boundary conditions.
Therefore $\Xi_1 = \ldots = \Xi_{p+1}$ and $\Xi_{n+1} = \ldots = \Xi_{n+p+1}$.
A set of B-Spline basis functions of degree two is illustrated in Fig.~\ref{fig:bs_basis}.

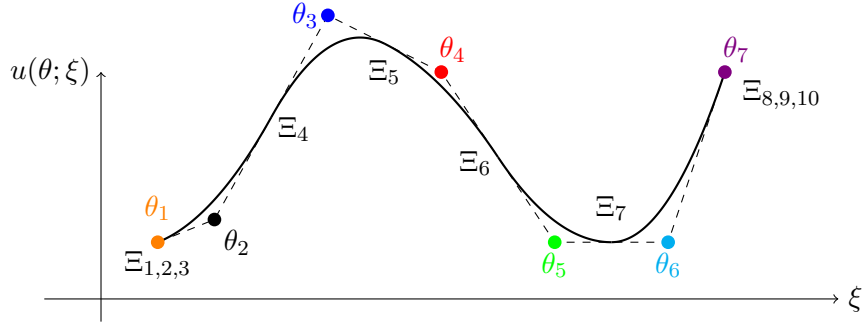
\begin{figure}[!ht]
    \begin{center}
    \begin{tikzpicture}[scale = 1.5]
    \draw[->] (-1,0) -- (6,0) node[right] {$\xi$};
    \draw[->] (-0.5,-0.2) -- (-0.5,2) node[left] {$u(\theta;\xi)$};

    \draw[dashed] (0,0.5) -- (0.5,0.7) -- (1.5,2.5) -- (2.5,2) -- (3.5,0.5) -- (4.5,0.5) -- (5,2);

    \draw[domain=0:1,smooth,variable=\x,thick] plot ({\x},{0.5*((1-\x)*(1-\x)) + 0.7*(2*\x-3/2*\x*\x) + 2.5*(0.5*\x*\x)});
    \node[orange] at (0,0.8) {$\theta_1$} ;
    \node[orange,circle,fill,scale=0.5] at (0,0.5) {};

    \draw[domain=1:2,smooth,variable=\x,thick] plot ({\x},{0.7*(0.5*(2-\x)*(2-\x)) + 2.5*(-0.5+\x-(\x-1)*(\x-1)) + 2*(0.5*(\x-1)*(\x-1))});
    \node[black] at (0.7,0.5) {$\theta_2$} ;
    \node[black,circle,fill,scale=0.5] at (0.5,0.7) {};

    \draw[domain=2:3,smooth,variable=\x,thick] plot ({\x},{2.5*(5/2-\x+0.5*(\x-2)*(\x-2)) + 2*(-0.5+(\x-1)-(\x-2)*(\x-2)) + 0.5*(0.5*(\x-2)*(\x-2))});
    \node[blue] at (1.3,2.5) {$\theta_3$} ;
    \node[blue,circle,fill,scale=0.5] at (1.5,2.5) {};

    \draw[domain=3:4,smooth,variable=\x,thick] plot ({\x},{2*(5/2-(\x-1)+0.5*(\x-3)*(\x-3)) + 0.5*(-0.5+(\x-2)-(\x-3)*(\x-3)) + 0.5*(0.5*(\x-3)*(\x-3))});
    \node[red] at (2.6,2.2) {$\theta_4$} ;
    \node[red,circle,fill,scale=0.5] at (2.5,2) {};

    \draw[domain=4:5,smooth,variable=\x,thick] plot ({\x},{0.5*(5/2-(\x-2)+0.5*(\x-4)*(\x-4)) + 0.5*(-16+10*(\x-1)-3/2*(\x-1)*(\x-1)) + 2*((\x-4)*(\x-4))});
    \node[green] at (3.5,0.3) {$\theta_5$} ;
    \node[green,circle,fill,scale=0.5] at (3.5,0.5) {};

    \node[cyan] at (4.5,0.3) {$\theta_6$} ;
    \node[cyan,circle,fill,scale=0.5] at (4.5,0.5) {};

    \node[violet] at (5.1,2.2) {$\theta_7$} ;
    \node[violet,circle,fill,scale=0.5] at (5,2) {};

    \node at (0,0.3) {$\Xi_{1,2,3}$} ;
    \node at (1.2,1.5) {$\Xi_4$} ;
    \node at (2,2) {$\Xi_5$} ;
    \node at (2.8,1.2) {$\Xi_6$} ;
    \node at (4,0.8) {$\Xi_7$} ;
    \node at (5.5,1.8) {$\Xi_{8,9,10}$} ;

    \end{tikzpicture}
    \end{center}
\caption{A quadratic B-Spline function.}
\label{fig:bs_curve}
\end{figure}

Finally, a one-dimensional B-Spline function is obtained by associating a set of degrees of freedom $(\theta_i)_{i=1,\cdots,n}$ to the basis functions:
\begin{equation}\label{eq:nurbscurve}
  u(\theta;\xi) = \sum_{i=1}^{n} B^p_i(\xi) \, \theta_i ,
\end{equation}
This representation is used in CAD to define the coordinate of a parametric curve $x(\theta;\xi)$. In the CAD context, $(\theta_i)_{i=1,\cdots,n}$ are the coordinates of control points in the \emph{physical domain} $\Omega$. Note that this equation defines a mapping $\mathcal{F}$ from the parametric domain $\widehat{\Omega}$ to the physical one $\Omega$. In the context of PDE solving, the degrees of freedom constitute an intuitive way to control locally the value of the function, as depicted in Fig~(\ref{fig:bs_curve}). On this plot, one can see the function $u(\theta;\xi)$ defined by the degrees of freedom $(\theta_i)_{i=1,\ldots,n}$ associated with the basis functions $(B^p_i)_{i=1,\ldots,n}$ shown previously, with the same colors.
Note that the polygon of degrees of freedom is tangent to the function $u$ at the knot parametric coordinates, which determine the support of each basis function. One can also notice that the function interpolates the degrees of freedom at both extremities and is tangent to the lines joining the two first and last degrees of freedom.

This representation can be extended to multi-dimensional variables by introducing a tensorial product of the univariate basis functions:
\begin{equation}\label{eq:nurbsmulti}
  u(\theta;\xi_1, \cdots, \xi_d) = \sum_{i_1=1}^{n_1} \cdots \sum_{i_d=1}^{n_d} B^{p_1}_{i_1}(\xi_1) \cdots B^{p_d}_{i_d}(\xi_d) \, \theta_{i_1, \cdots, i_d} .
\end{equation}

The time and space derivatives of the solution can be simply obtained by analytical differentiation of the function definition. Thus, the derivative of the solution with respect to the variable $\xi_j$ writes:
\begin{equation}\label{eq:nurbsderiv}
  \partial_{\xi_j} u(\theta;\xi_1, \cdots, \xi_j, \cdots, \xi_d) = \sum_{i_1=1}^{n_1} \cdots \sum_{i_d=1}^{n_d} B^{p_1}_{i_1}(\xi_1) \cdots [B^{p_j}_{i_j}]'(\xi_j) \cdots B^{p_d}_{i_d}(\xi_d) \, \theta_{i_1, \cdots, i_d} ,
\end{equation}
where $B'$ denotes for the derivative of the univariate B-Spline basis function. Its computation can be efficiently implemented using a dedicated algorithm, based on the fact that the derivative of a B-Spline can be expressed as a B-Spline of lower degree~\cite{DeBoor_78}. Note that evaluating the derivatives in the Cartesian frame can be more complex, due to the possible composition with a geometrical mapping. This point will be detailed in the presentation of the last testcase.

The differentiation of the loss function with respect to the degrees of freedom $\theta$ can also be achieved analytically. It relies on the derivatives of the solution, which simply writes:
\begin{equation}\label{eq:nurbsderivtheta}
  \partial_{\theta{i_1, \cdots, i_d}} u(\theta;\xi_1, \cdots, \xi_d) =  B^{p_1}_{i_1}(\xi_1) \cdots B^{p_d}_{i_d}(\xi_d) .
\end{equation}

As explained above, the imposition of the initial / boundary values and derivatives is straightforward when clamped knot vectors are used. Indeed, in that case, the extremity of the B-Spline is only dependent on the extreme degrees of freedom, as shown in Fig.~(\ref{fig:bs_curve}). For instance, in the case of a two-parameter B-Spline surface, the value of the B-Spline along the curve $(\xi_1=\Xi_1,\xi_2)$ is only defined by the $n_2$ degrees of freedom $(\theta_{1,i_2})_{i_2 = 1,\ldots,n_2}$. The value of the solution on this boundary can therefore be enforced independently from the training precedure by using a least-square projection:
\begin{equation}\label{eq:nurbsapprox}
  (\theta_{1,i_2})_{i_2 = 1,\ldots,n_2} = \text{ArgMin}  \int_{\partial \Omega_L} \left [ u(\theta;\Xi_1, \xi_2) - g^D(\xi_2) \right ]^2 d\xi_2 ,
\end{equation}
where $\partial \Omega_L$ is the boundary defined by $\xi_1 = \Xi_1$. This problem is quadratic in $\theta$ and yields a linear system~\cite{DeBoor_78}. Note that the boundary solution is exact in this example if the Dirichlet condition $g^D(\xi_2)$ is polynomial of degree $p_2$ or less. For Neumann conditions, a similar approach can be adopted, involving the first degree of freedom inside the domain to express the boundary slope, as illustrated in Fig.~(\ref{fig:bs_curve}).

\section{Properties of the two representations}
\label{sec:prop}

In this section, we compare the properties of the two representations, in the perspective of the physics-informed learning procedure detailed above.

\paragraph{Support}

A first major difference between the two representations is the domain of influence of the degrees of freedom. For B-Spline bases, each function is non-zero over $p+1$ knot intervals, as illustrated by Fig.~\ref{fig:bs_basis}. In contrast, functions generated by MLPs have infinite support since they are compositions of $tanh$ activation functions. Consequently, the influence of each individual parameter $\theta_i$ is located in space and time for B-Spline representations, and not for MLP networks. Thus, some evaluations can be avoided for B-Splines when the collocation points are located outside the support of a given basis function.

\paragraph{Smoothness}

Since a $C^\infty$ function is usually employed as activation, e.g. $tanh$, the output of MLP networks are also $C^\infty$. B-Spline functions are $C^\infty$, except at the knots abscissae where a $C^{p-1}$ continuity is obtained by construction, if the knots are distinct~\cite{DeBoor_78}.

\paragraph{Linear independence}

The B-Spline functions form a polynomial basis, then they are linearly independent and generate functions depending linearly on the degrees of freedom. The functions generated by the hidden layers of MLP networks are obviously non-linear with respect to $\theta$ and not independent in general.

\paragraph{Convexity}

For some problems, it is easy to prove the convexity of the loss function using a B-Spline representation. For instance, if one considers a one-dimensional Poisson equation $\partial_{xx} u = f$ associated with Dirichlet boundary conditions, the loss function simply writes $\frac{1}{N}\sum_k | \mathbf{B}''(\xi_k)^T \, \theta - f(\xi_k) |^2$ according to Eq.~(\ref{eq:pdeloss}) and~(\ref{eq:nurbscurve}), with $\mathbf{B}''(\xi_k)$ a vector collecting the second-order derivatives of the basis functions. Then, the Hessian matrix of the loss function is $\frac{2}{N} \sum_k \mathbf{B}''(\xi_k) \, \mathbf{B}''(\xi_k)^T$, which is positive-definite, for an adequate choice of collocation points.
Therefore, the problem admits a unique solution.
On the contrary, the loss functions obtained using MLP networks are globaly non-convex. Indeed, swapping the weights and biases of a network can result in the exact same output, yielding several minima.

\paragraph{Approximation capabilities}

It has been shown that MLP networks have the ability to represent a large variety of functions. In particular, a neural network with a single hidden layer and a non-polynomial activation function can approximate any continuous function to any desired accuracy~\cite{hornikMultilayerFeedforwardNetworks1989}. Besides, B-Splines are known to have noticeable  approximation properties. In particular, the error for approximating a smooth function decays at the optimal rate of $O(h^{p+1})$ in the $L^{\infty }$ and $L^2$ norms~\cite{DeBoor_78}, where $h$ is the knot interval.

\paragraph{Computation of derivatives}

As already explained, the computation of the derivatives for neural networks is achieved using AD techniques, whereas analytic expressions of derivatives can be used for B-Splines. Obviously, the recourse to AD facilitates the implementation, however the evaluation of analytic derivatives is usually far more efficient in terms of computational time and memory use.

\paragraph{Account for geometry}

The neural networks consider as input the space coordinates, usually expressed in Cartesian frame for the sake of simplicity. Nevertheless, other frames, e.g. cylindrical or spherical, can be used as well if more convenient. The situation is different for B-Splines, because the basis functions are defined in the parametric domain, which should be mapped to the physical domain by fitting. This mapping makes the expression of the derivatives in the Cartesian frame more complex, because the geometric mapping has to be differentiated as well. Nevertheless, this allows to match the boundary of the physical domain to B-Spline extremities, thus facilitating the imposition of the boundary conditions. This point will be especially illustrated in the last testcase.

\paragraph{Imposition of boundary conditions}

In the PINNs methodology, the imposition of the initial and boundary conditions is achieved by penalizing the residual term in the loss function, yielding a multi-objective training problem. It has been shown that the antagonism between the different loss terms is a cause of ill-conditioning~\cite{krishnapriyan2021characterizing,wangExpertGuideTraining2023}.
In contrast, the properties of B-Splines at the domain extremities, by the use of clamped knot vectors, allow to impose intial and boundary conditions by projection, independently from the training phase. This is especially important, because it reduces the loss function to the PDE residual only and, consequently, facilitates the training. Moreover, the degrees of freedom corresponding to the B-Spline extremities can be removed from the set of trained parameters, yielding a problem of reduced size.

\paragraph{Sampling}

For MLP networks, collocation points are usually chosen randomly, according to a given distribution that promotes space filling, e.g. Latin Hypercube Sampling (LHS). Re-sampling or mini-batches are commonly employed during the training phase, emphasizing the stochastic nature of the training. Regarding B-Spline representations, a similar strategy could obviously be adopted. However, the specific structure of B-Spline bases allows an alternate choice based on Greville abscissae, which are the mean locations of $p$ consecutive knots in the knot vector $g_k = \frac{1}{p}\sum_{i=1}^p \Xi_{k+i}$. It has been shown that this set of points, whose cardinal is equal to the number of functions, allow a stable B-Spline interpolation and yields orders of convergence $p$ and $p-1$ for even and odd degrees respectively~\cite{JOHNSON200563}. The resulting approach relies only on $n$ collocation points per direction and is in fact similar to the isogeometric collocation method~\cite{Reali2015}.

\paragraph{Initialization}

The inital guess of the solution can be critical for the training, since descent algorithms are essentially local optimizers. The initialization of network parameters is usually achieved randomly with a dedicated distribution, e.g. Xavier uniform initialization~\cite{pmlr-v9-glorot10a}, because it is difficult to provide a guess for the network parameters. As consequence, the robustness of the resolution has to be assessed by performing several runs because different solutions may be found in practice depending on the initialization. On the contrary, the B-Spline solution can be initialized by intuiting a solution guess, since there is a direct relationship between the B-Spline degrees of freedom and the solution, as illustrated in Fig.~(\ref{fig:bs_curve}).

\medskip

Following this review of the properties of the two chosen representations, a set of computational testcases are presented and analysed in the next section.

\section{Computational Tests}
\label{sec:result}

\subsection{Numerical setup}

A set of problems with closed-form solution is now solved using both representations. Some numerical parameters common to all computations are detailed below. They have been selected according to the previous discussions and similar studies from the literature.

\paragraph{Neural representation}
For MLP networks, one considers a set of different architectures, ranging from one to five hidden layers, composed of 8, 16, 32 and 64 neurons each. The activation function is always $tanh$. The sampling of the collocation points is achieved using a LHS technique and remains fixed during the whole training. The number of points is 2000 in the domain, plus 100 on each boundary. All boundary conditions are implemented as penalty terms in the loss function, with unitary weights. The network parameters are initialized with a Xavier uniform distribution. For each test, a set of 5 runs is performed to assess the robustness with respect to the initialization and sampling. The training is organized in two distinct phases: first, ADAM optimizer is used during 2000 iterations (fixed learning rate $10^{-3}$), followed by BFGS optimizer until full convergence of the loss function.

\paragraph{Spline representation}
For B-Splines, one considers basis degrees ranging from 2 to 6, with an increasing number of degrees of freedom in each direction, depending on the case. A uniform knot vector is employed, with clamped extremities. All boundary conditions are taken into account by projection, except for the last testcase. The Greville abscissae are chosen as collocation points. The B-Splines are initialized with a constant value depending on the problem. The training is here performed in a single phase, using directly BFGS optimizer, until full convergence.

\paragraph{Error estimate}

The accuracy of the resolution is measured by computing the Root Mean Square Error (RMSE) and the Maximum Absolute Error (MaxAE).
RMSE is estimated by computing the difference between the analytical solution $u^\star$ and the numerical one $u$, for a uniform grid including $N_\epsilon = 512$ points per direction:
\begin{equation}
    \epsilon_2(\theta) = \sqrt{\frac{1}{N_\epsilon^2}\sum_{i=1}^{N_\epsilon} \sum_{j=1}^{N_\epsilon} \left | u(\theta;x_i, t_j) - u^\star(x_i, t_j) \right |^2}
\end{equation}
The same grid is employed for the two representations, ensuring a fair comparison.
MaxAE is also estimated using the same grid:
\begin{equation}
    \epsilon_\infty(\theta) = \max_{i,j=1,\cdots,N_\epsilon} \left | u(\theta;x_i, t_j) - u^\star(x_i, t_j) \right |
\end{equation}

For neural representations, since a set of independant runs are performed for each test, the average, maximum and minimum values obtained are provided. Moreover, the absolute error fields $| u(\theta) - u^\star|$ obtained using the two representations are compared in all testcases, for configurations counting about 1200 degrees of freedom.

\paragraph{Computational time}

The comparison of the computational times required to reach a target accuracy would be obviously of great interest. However, runtime depends on several factors, like hardware characteritics, implementation, code language, etc. that make the comparison tedious and usually impossible to reproduce. Therefore, we do not provide computational times in the present work and prefer focus on error estimates.

\clearpage

\subsection{Poisson equation}

\begin{figure}[!ht]
    \centering
    \includegraphics[width=0.5\textwidth]{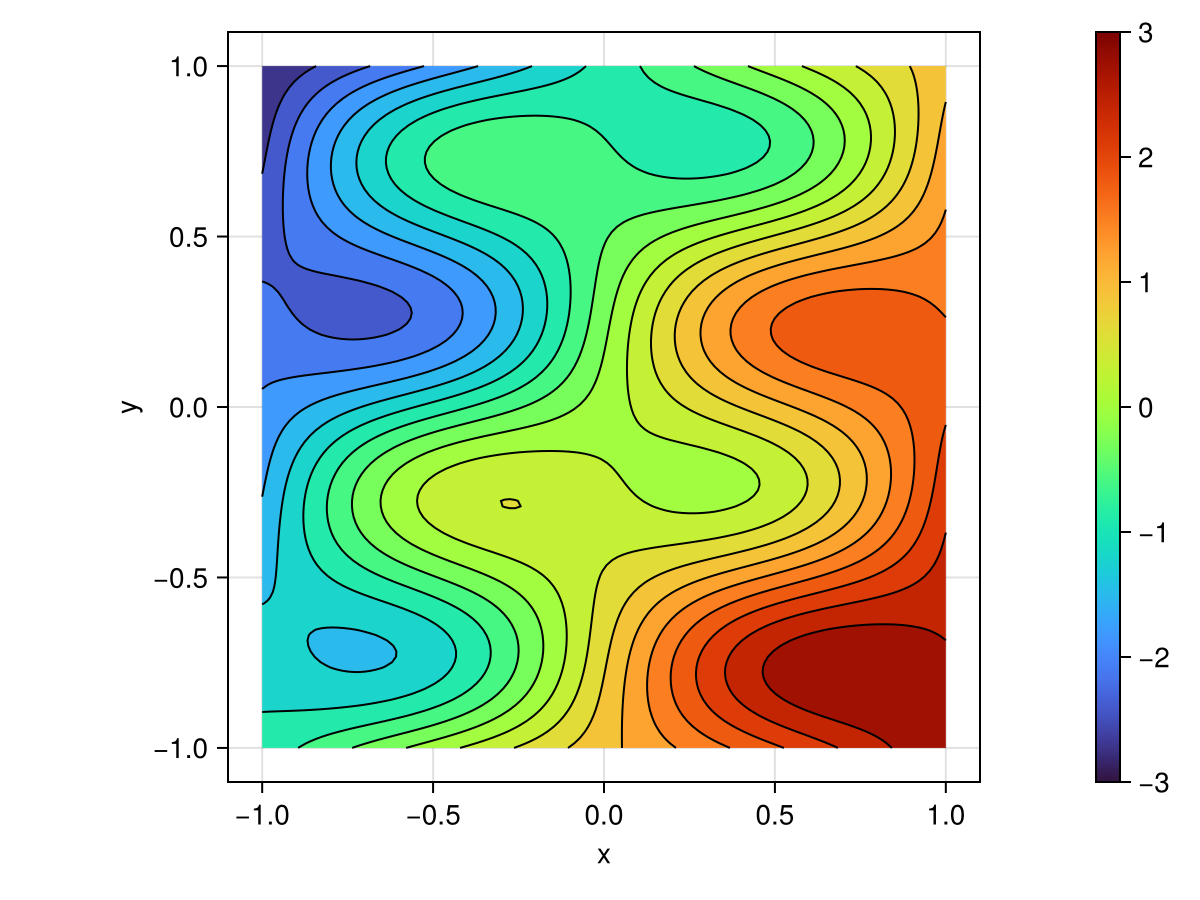}
     \caption{Analytical solution for Poisson equation.}
    \label{fig:poisson_exact}
\end{figure}

\begin{figure}[!ht]
  \begin{subfigure}[b]{0.49\linewidth}
    \centering
    \includegraphics[width=\textwidth]{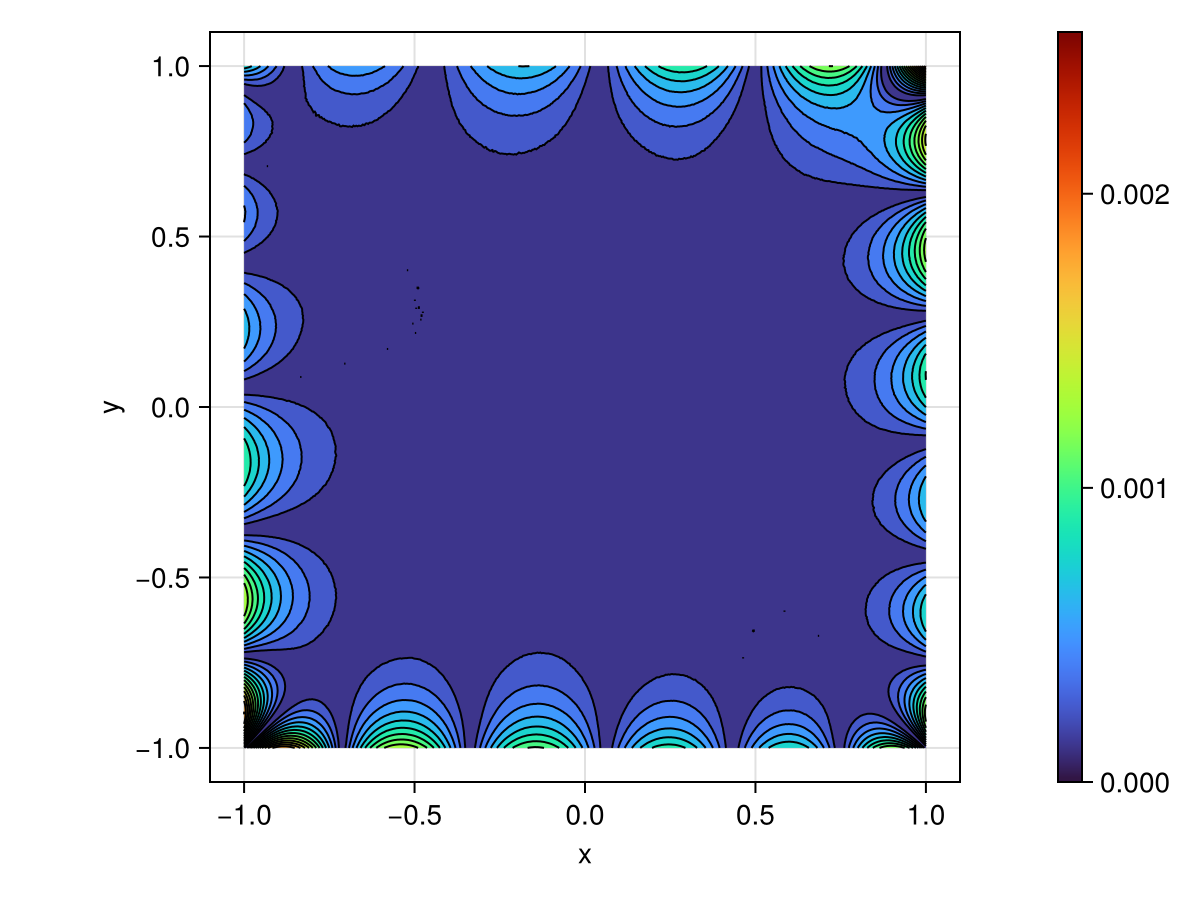}
     \caption{Neural network}
    \label{fig:poisson2d_error_nn}
  \end{subfigure}
  \begin{subfigure}[b]{0.49\linewidth}
    \centering
    \includegraphics[width=\textwidth]{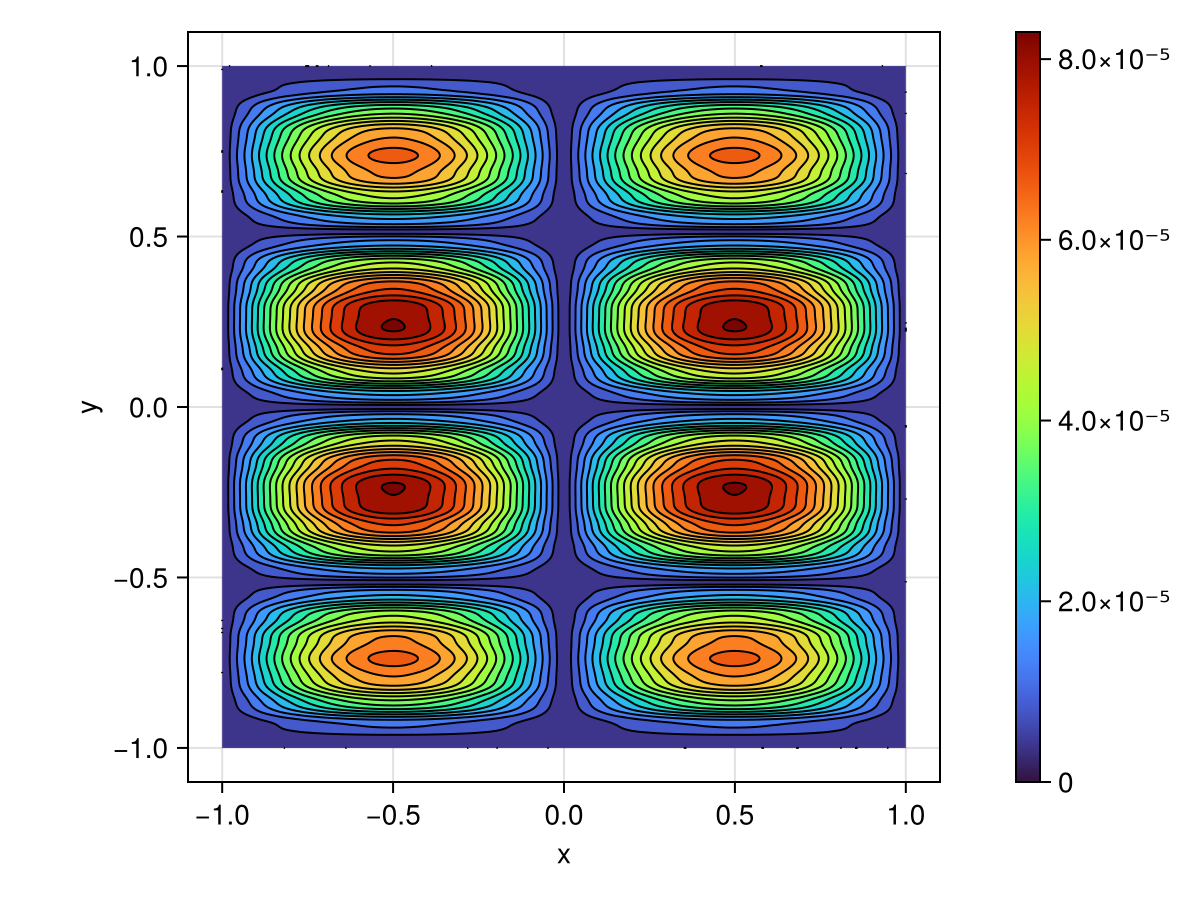}
    \caption{B-Spline basis}
    \label{fig:poisson2d_error_sp}
  \end{subfigure}
  \caption{Error fields for Poisson equation (about 1200 parameters).}
  \label{fig:poisson_error}
\end{figure}

\begin{figure}[!ht]
  \begin{subfigure}[b]{0.49\linewidth}
    \centering
    \includegraphics[width=\textwidth]{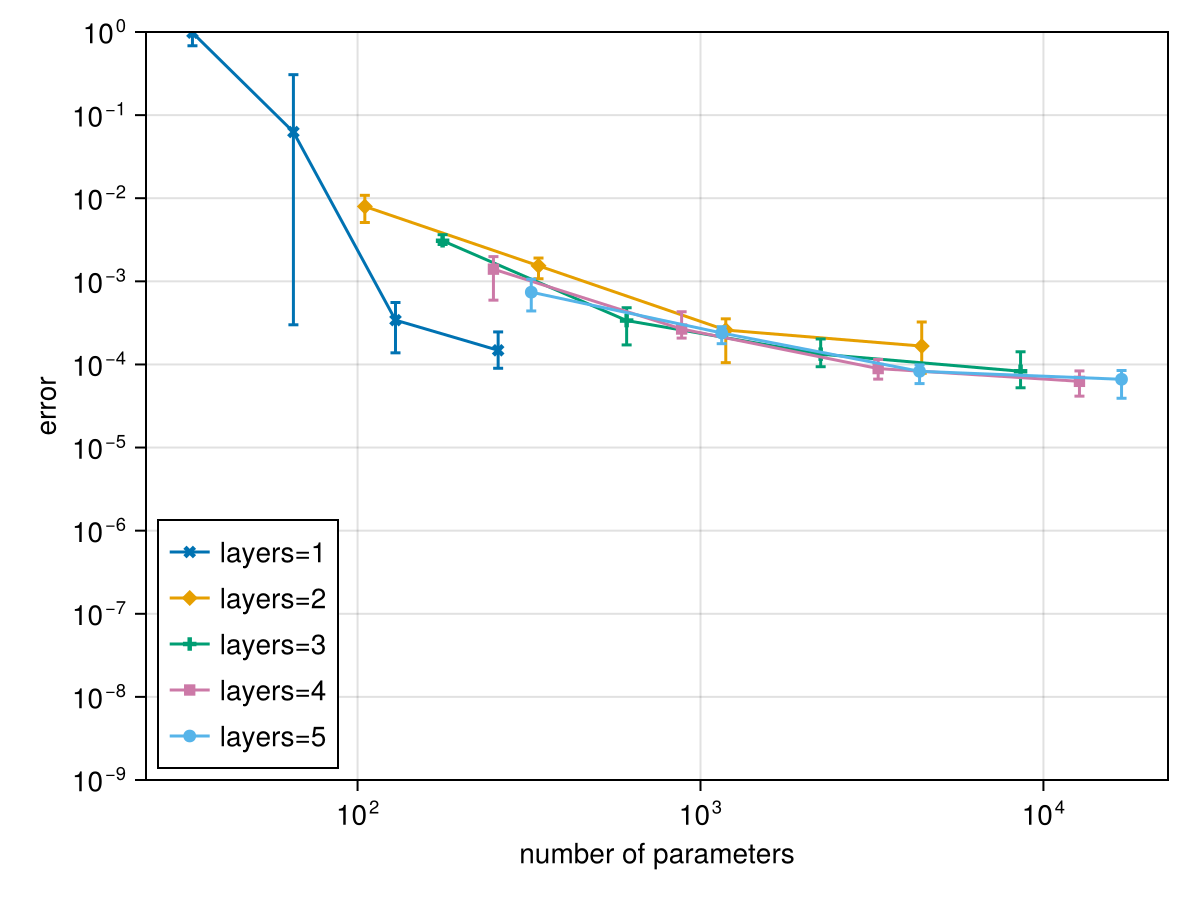}
     \caption{Neural network}
    \label{fig:poisson_nn_l2}
  \end{subfigure}
  \begin{subfigure}[b]{0.49\linewidth}
    \centering
    \includegraphics[width=\textwidth]{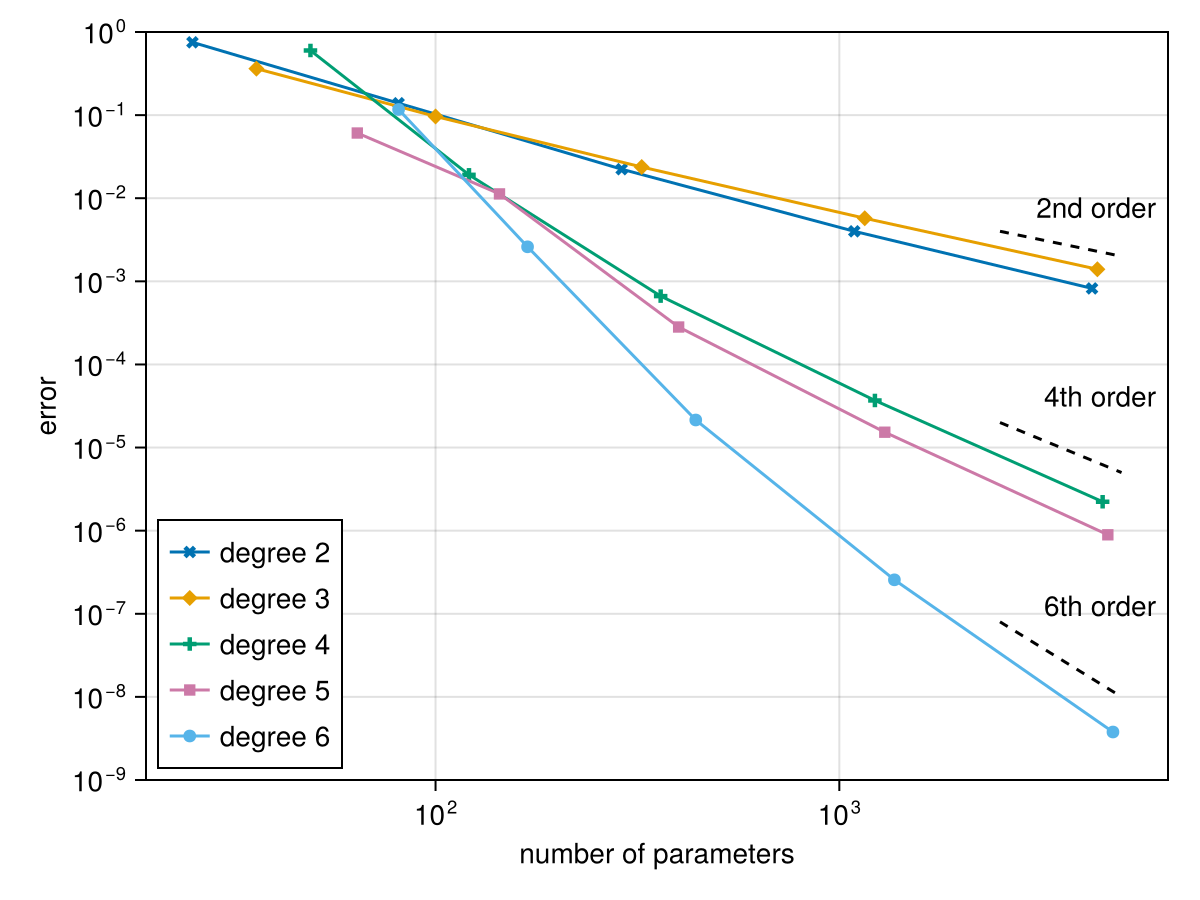}
    \caption{B-Spline basis}
    \label{fig:poisson_sp_l2}
  \end{subfigure}
  \caption{$\epsilon_2$ norm of error w.r.t. the number of degrees of freedom for Poisson equation.}
  \label{fig:poisson_l2}
\end{figure}

The first problem studied is a two-dimensional linear elliptic problem:
\begin{equation}
    \partial_{xx} u + \partial_{yy} u = 5 \pi^2 \sin(\pi x) \sin(2 \pi y) \quad \forall (x,y) \in ]-1,1[^2 ,
\end{equation}
with the boundary conditions:
\begin{align}
    u(x,-1) &= 2x+1 \\
    u(x,1) &= 2x-1 \\
    u(-1,y) &= -y-2 \\
    u(1,y) &= -y+2 .
\end{align}
This problem admits as closed-form solution $u(x,y) = 2x -y + \sin(\pi x) \sin(2 \pi y)$.
For B-Splines, the solution is initialized with a zero field, except the degrees of freedom located at the boundaries, which are obtained by projection of the analytical solution. B-Spline representations are tested using 4, 8, 16, 32 and 64 degrees of freedom in each direction.
The analytical solution can be seen in Fig.~(\ref{fig:poisson_exact}) and the variation in error as a function of the number of degrees of freedom is depicted in Fig.~(\ref{fig:poisson_l2}), with the same log scale for the two representations. Numerical values of error are provided in Tabs.~(\ref{tab:poisson_nn}-\ref{tab:poisson_sp}). The absolute error fields are shown in Fig.~(\ref{fig:poisson_error}).

\begin{table}
    \centering
    \begin{tabular}{|c|c|c||c|c|c|c|c|c| }
        \hline
        Layers & Neurons & d.o.f. & \multicolumn{3}{|c|}{$\epsilon_2$} & \multicolumn{3}{|c|}{$\epsilon_\infty$} \\
        \hline
         & & & min & average & max & min & average & max \\
        \hline
        \hline
        1 & 8 & 33 & 6.83e-01 & 9.81e-01 & 1.26e+00 & 3.15e+00 & 4.06e+00 & 4.84e+00 \\
        1 & 16 & 65 & 3.00e-04 & 6.26e-02 & 3.07e-01 & 1.03e-03 & 5.21e-01 & 2.58e+00 \\
        1 & 32 & 129 & 1.38e-04 & 3.41e-04 & 5.56e-04 & 6.05e-04 & 1.59e-03 & 2.85e-03 \\
        1 & 64 & 257 & 8.99e-05 & 1.48e-04 & 2.46e-04 & 5.27e-04 & 9.36e-04 & 1.30e-03 \\
        \hline
        \hline
        2 & 8 & 105 & 5.10e-03 & 7.96e-03 & 1.08e-02 & 3.29e-02 & 5.55e-02 & 7.36e-02 \\
        2 & 16 & 337 & 1.07e-03 & 1.54e-03 & 1.91e-03 & 8.81e-03 & 1.50e-02 & 1.74e-02 \\
        2 & 32 & 1185 & 1.05e-04 & 2.59e-04 & 3.54e-04 & 7.87e-04 & 2.50e-03 & 3.44e-03 \\
        2 & 64 & 4417 & 7.83e-05 & 1.66e-04 & 3.24e-04 & 8.84e-04 & 1.53e-03 & 2.77e-03 \\
        \hline
        \hline
        3 & 8 & 177 & 2.78e-03 & 3.08e-03 & 3.64e-03 & 2.12e-02 & 3.02e-02 & 3.36e-02 \\
        3 & 16 & 609 & 1.72e-04 & 3.37e-04 & 4.81e-04 & 1.92e-03 & 3.32e-03 & 4.14e-03 \\
        3 & 32 & 2241 & 9.39e-05 & 1.33e-04 & 2.03e-04 & 9.61e-04 & 1.14e-03 & 1.34e-03 \\
        3 & 64 & 8577 & 5.24e-05 & 8.28e-05 & 1.42e-04 & 4.86e-04 & 9.57e-04 & 1.42e-03 \\
        \hline
        \hline
        4 & 8 & 249 & 5.94e-04 & 1.40e-03 & 1.98e-03 & 4.09e-03 & 1.19e-02 & 1.83e-02 \\
        4 & 16 & 881 & 2.08e-04 & 2.67e-04 & 4.31e-04 & 1.93e-03 & 2.88e-03 & 5.06e-03 \\
        4 & 32 & 3297 & 6.66e-05 & 8.92e-05 & 1.15e-04 & 7.02e-04 & 9.32e-04 & 1.24e-03 \\
        4 & 64 & 12737 & 4.16e-05 & 6.26e-05 & 8.36e-05 & 5.26e-04 & 7.23e-04 & 1.09e-03 \\
        \hline
        \hline
        5 & 8 & 321 & 4.40e-04 & 7.40e-04 & 1.06e-03 & 4.48e-03 & 8.87e-03 & 1.29e-02 \\
        5 & 16 & 1153 & 1.78e-04 & 2.39e-04 & 2.85e-04 & 1.81e-03 & 2.36e-03 & 2.92e-03 \\
        5 & 32 & 4353 & 5.88e-05 & 8.26e-05 & 9.89e-05 & 5.85e-04 & 7.78e-04 & 1.00e-03 \\
        5 & 64 & 16897 & 3.92e-05 & 6.63e-05 & 8.46e-05 & 3.32e-04 & 6.80e-04 & 9.97e-04 \\
        \hline
    \end{tabular}
    \caption{Errors for Poisson equation for MLP networks.}
    \label{tab:poisson_nn}
\end{table}

\begin{table}
    \centering
    \begin{tabular}{|c|c||c|c|}
        \hline
        Degree & d.o.f. & $\epsilon_2$ & $\epsilon_\infty$ \\
        \hline
        \hline
        2 & 25 & 7.51e-01 & 1.83e+00 \\
        2 & 81 & 1.39e-01 & 3.94e-01 \\
        2 & 289 & 2.23e-02 & 6.41e-02 \\
        2 & 1089 & 4.00e-03 & 9.13e-03 \\
        2 & 4225 & 8.22e-04 & 1.73e-03 \\
        \hline
        \hline
        3 & 36 & 3.62e-01 & 9.08e-01 \\
        3 & 100 & 9.63e-02 & 2.55e-01 \\
        3 & 324 & 2.37e-02 & 5.01e-02 \\
        3 & 1156 & 5.74e-03 & 1.16e-02 \\
        3 & 4356 & 1.40e-03 & 2.82e-03 \\
        \hline
        \hline
        4 & 49 & 5.99e-01 & 1.52e+00 \\
        4 & 121 & 1.91e-02 & 4.36e-02 \\
        4 & 361 & 6.64e-04 & 1.64e-03 \\
        4 & 1225 & 3.69e-05 & 8.35e-05 \\
        4 & 4489 & 2.23e-06 & 4.77e-06 \\
        \hline
        \hline
        5 & 64 & 6.09e-02 & 1.69e-01 \\
        5 & 144 & 1.13e-02 & 3.31e-02 \\
        5 & 400 & 2.82e-04 & 7.13e-04 \\
        5 & 1296 & 1.53e-05 & 3.21e-05 \\
        5 & 4624 & 8.92e-07 & 1.83e-06 \\
        \hline
        \hline
        6 & 81 & 1.17e-01 & 2.78e-01 \\
        6 & 169 & 2.60e-03 & 8.13e-03 \\
        6 & 441 & 2.15e-05 & 5.59e-05 \\
        6 & 1369 & 2.56e-07 & 6.17e-07 \\
        6 & 4761 & 3.78e-09 & 8.95e-09 \\
        \hline
    \end{tabular}
    \caption{Errors for Poisson equation for B-Splines.}
    \label{tab:poisson_sp}
\end{table}

As can be seen, the convergence of the error obtained using B-Spline representations corresponds perfectly to the expected theoretical results. The error tends to zero as the number of parameters increases and convergence rates of value $p$ and $p-1$ are observed for even and odd degrees respectively. Results obtained using neural networks are very different. Whatever the network architecture, the error trajectories follow a similar trend, with a fast reduction of error for small networks and a low one for larger networks. It should be noticed that making the network deeper or larger yields similar effects. Both methods can reach the error level $10^{-4}$ with a few parameters (about $500$), but neural networks fail in reducing the error below this threshold, contrary to B-Splines. As shown by the absolute error maps in Fig.~(\ref{fig:poisson_error}), the error for the network is localized at the boundaries of the domain, which highlights an unbalanced loss formulation. For B-Splines, the tensorial nature of the representation can be obviously observed.

\clearpage

\subsection{Burgers equation with viscosity}

\begin{figure}[!ht]
    \centering
    \includegraphics[width=0.5\textwidth]{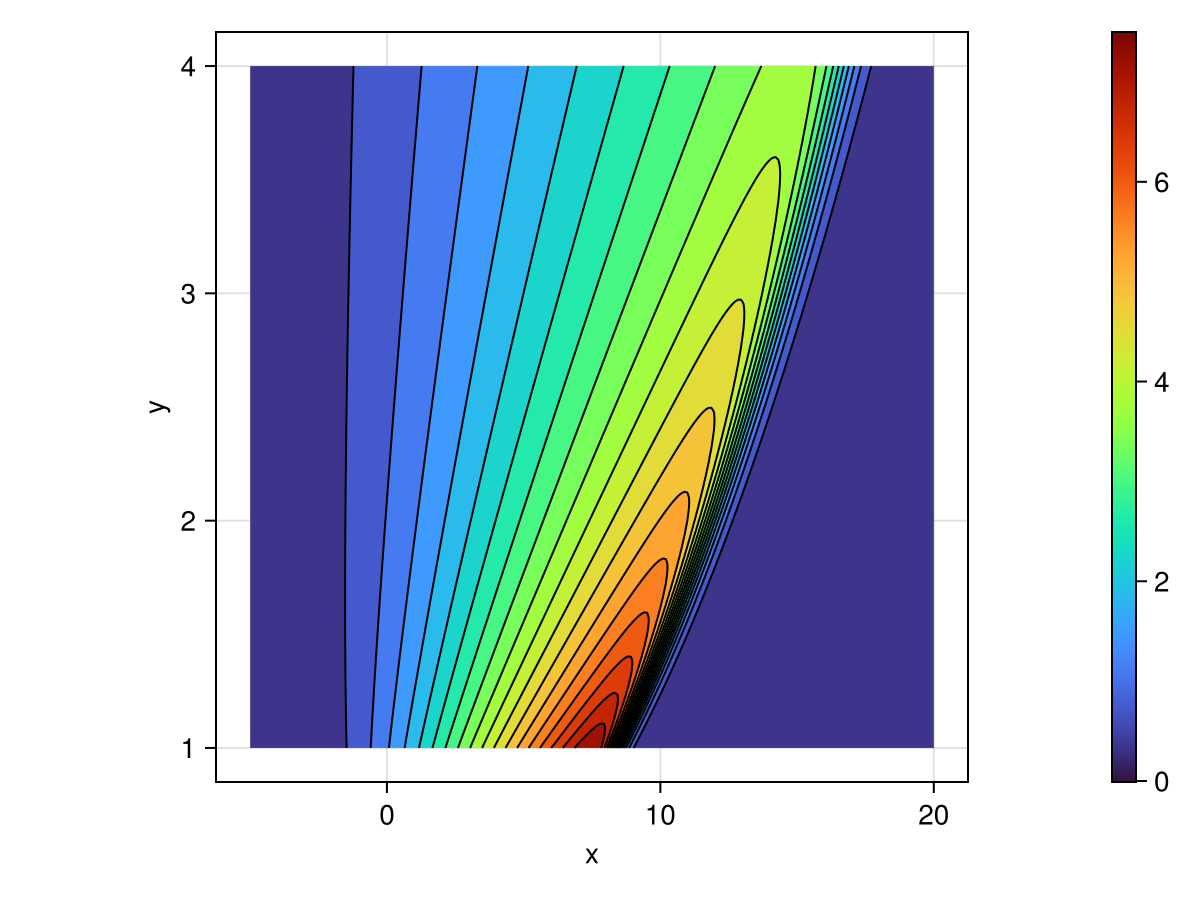}
     \caption{Analytical solution for Burgers equation with viscosity.}
    \label{fig:burgers_exact}
\end{figure}

As second problem, one considers the Burgers equation with a viscous term, which allows to introduce a non-linear hyperbolic term, whose influence is balanced by adjusting the viscosity coefficient $\nu$:
\begin{equation}
    \partial _t u + u \, \partial_x u -\nu \partial_{xx} u = 0 \qquad \forall (x,t) \in ]-5,20[ \times ]1,4] .
\end{equation}
This equation admits a closed-form solution for $t>0$, for an initial state equal to a Dirac distribution $2\nu R \delta(x)$ at $t=0$~\cite{wave_11}:
\begin{equation}
    u^\star(x,t) = \sqrt{\frac{\nu}{\pi t}} \; \frac{(e^R -1) \, e^{-\frac{x^2}{4 \nu t}}}{1+\frac{1}{2}(e^R -1) \, \text{erfc}(\frac{x}{4 \nu t})} ,
\end{equation}
with erfc() the complement of the cumulative distribution function.
The problem considered numerically is solved for the time interval $]1,4]$ with $\nu=1$ and $R=20$. The initial condition at $t=1$ and the boundary conditions at $x=-5$ and $x=20$ are set using the analytical solution, by projection for B-Splines and penalization for neural networks. The remaining B-Spline parameters are initialized with zero value. B-Spline representations are tested using 8, 16, 32, 64 and 128 parameters in each direction. The analytical solution is depicted in Fig.~(\ref{fig:burgers_exact}). It is characterized by a smooth shock traveling from left to right, with a progressive dissipation as time growth.

\begin{figure}[!ht]
  \begin{subfigure}[b]{0.49\linewidth}
    \centering
    \includegraphics[width=\textwidth]{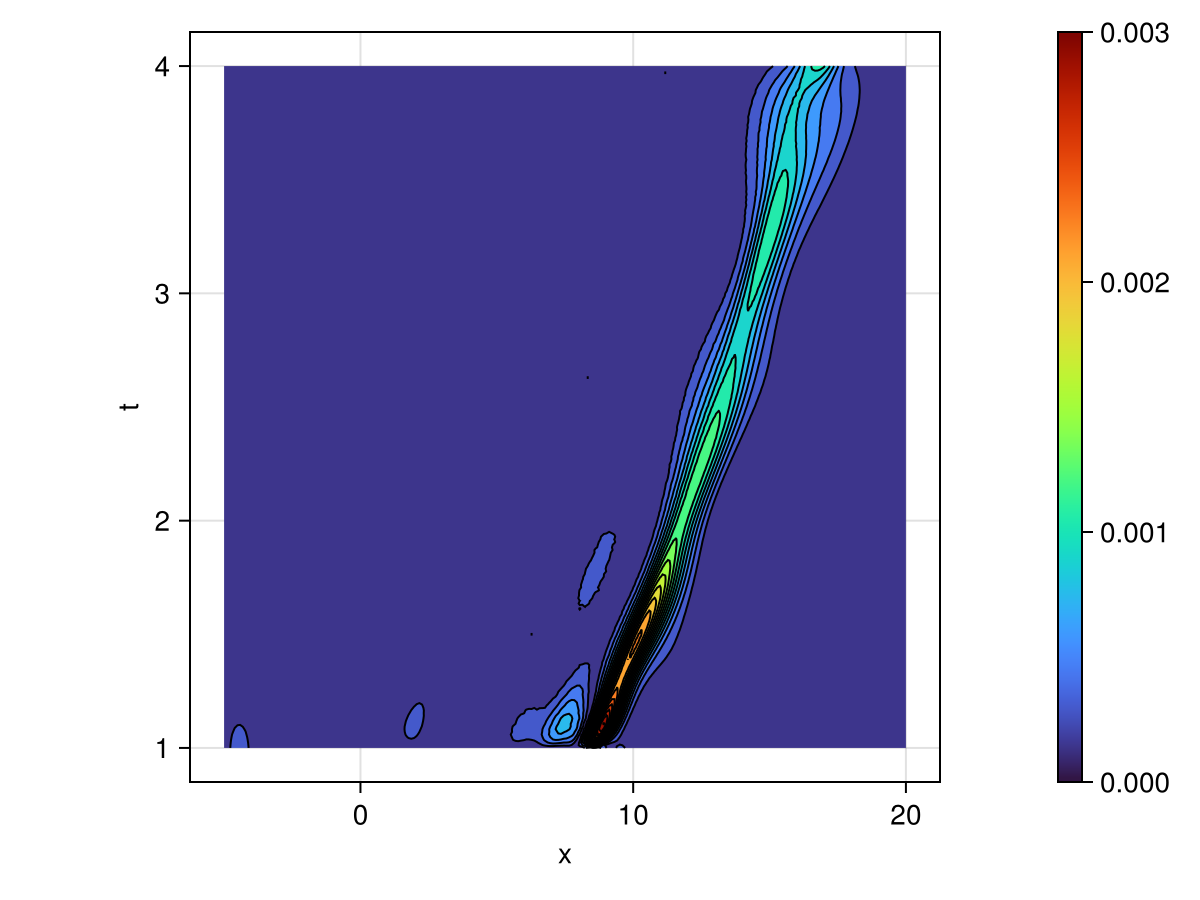}
     \caption{Neural network}
    \label{fig:burgers_error_nn}
  \end{subfigure}
  \begin{subfigure}[b]{0.49\linewidth}
    \centering
    \includegraphics[width=\textwidth]{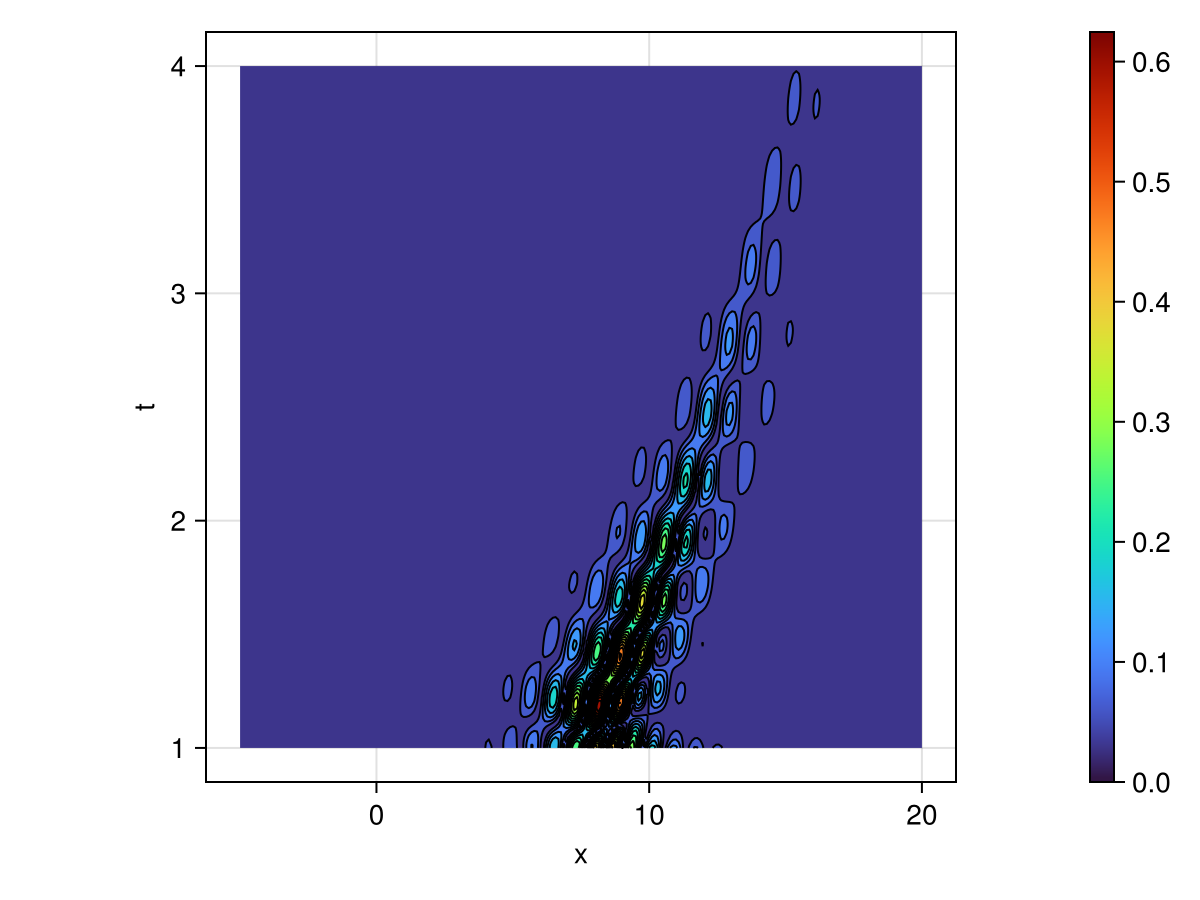}
    \caption{B-Spline basis}
    \label{fig:burgers_error_sp}
  \end{subfigure}
  \caption{Error fields for Burgers equation with viscosity (about 1200 parameters).}
  \label{fig:burgers_error}
\end{figure}

\begin{figure}[!ht]
  \begin{subfigure}[b]{0.49\linewidth}
    \centering
    \includegraphics[width=\textwidth]{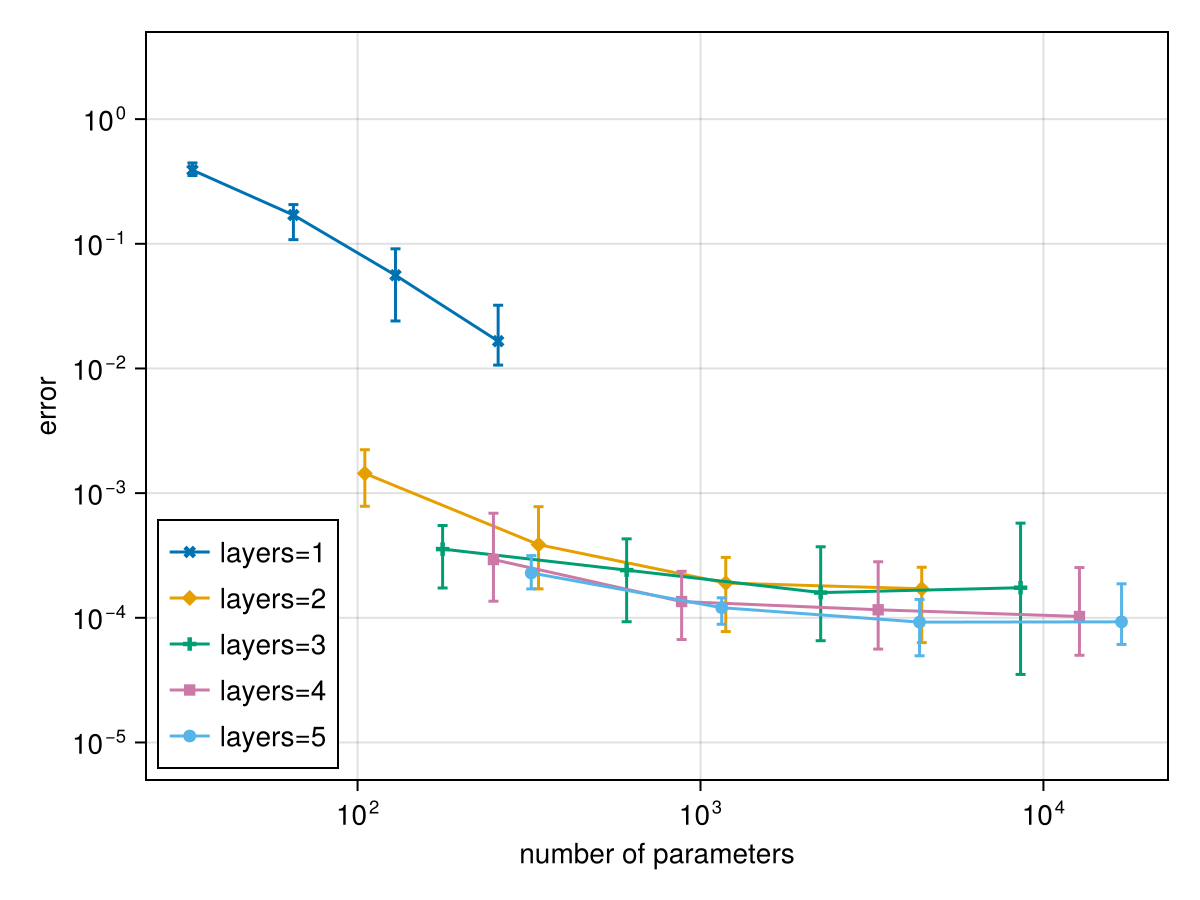}
     \caption{Neural network}
    \label{fig:burgers_nn_l2}
  \end{subfigure}
  \begin{subfigure}[b]{0.49\linewidth}
    \centering
    \includegraphics[width=\textwidth]{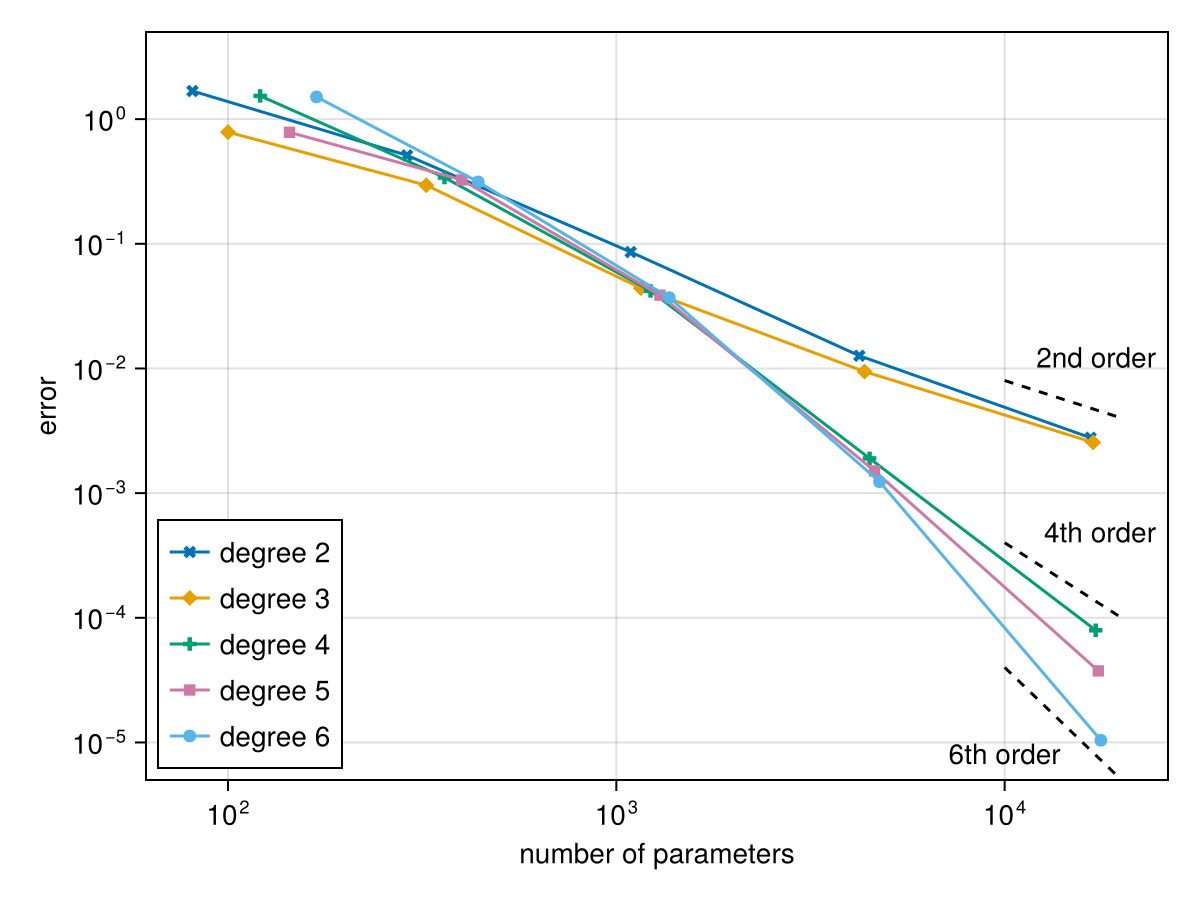}
    \caption{Spline basis}
    \label{fig:burgers_sp_l2}
  \end{subfigure}
  \caption{$\epsilon_2$ norm of error w.r.t. the number of parameters for Burgers equation with viscosity.}
  \label{fig:burgers_l2}
\end{figure}

The variation in error as a function of the number of parameters is depicted in Fig.~(\ref{fig:burgers_l2}) and provided in Tabs.~(\ref{tab:burgers_nn}-\ref{tab:burgers_sp}). As can be seen, the presence of the non-linear hyperbolic term makes the problem more difficult. The method using B-Splines exhibits a high-level error when the number of parameters is low, whatever the degree employed, but succeeds in reducing the error when the bases are refined enough. The asymptotic convergence rates are approximately recovered. Error maps are presented in Fig.~(\ref{fig:burgers_error}). Actually, the error is found located around the shock (Gibbs phenomenon), which cannot be solved accurately as long as there are not enough degrees of freedom in its vicinity. The situation is quite different using neural networks, that can obtain a good approximation of the solution with a low number of parameters, except single-layer networks. However, as for the previous problem, there seems to be a plateau at error level $10^{-4}$, that cannot be exceeded. One can also notice that the dispersion of results obtained for each test is larger than the one observed for the previous problem.

\begin{table}
    \centering
    \begin{tabular}{|c|c|c||c|c|c|c|c|c| }
        \hline
        Layers & Neurons & d.o.f. & \multicolumn{3}{|c|}{$\epsilon_2$} & \multicolumn{3}{|c|}{$\epsilon_\infty$} \\
        \hline
         & & & min & average & max & min & average & max \\
        \hline
        \hline
        1 & 8 & 33 & 3.53e-01 & 3.90e-01 & 4.45e-01 & 2.68e+00 & 2.96e+00 & 3.22e+00 \\
        1 & 16 & 65 & 1.08e-01 & 1.70e-01 & 2.06e-01 & 6.73e-01 & 1.11e+00 & 1.37e+00 \\
        1 & 32 & 129 & 2.40e-02 & 5.59e-02 & 9.11e-02 & 1.82e-01 & 4.16e-01 & 6.85e-01 \\
        1 & 64 & 257 & 1.06e-02 & 1.66e-02 & 3.22e-02 & 8.15e-02 & 1.51e-01 & 3.48e-01 \\
        \hline
        \hline
        2 & 8 & 105 & 7.85e-04 & 1.44e-03 & 2.23e-03 & 6.23e-03 & 1.03e-02 & 1.54e-02 \\
        2 & 16 & 337 & 1.71e-04 & 3.86e-04 & 7.79e-04 & 1.65e-03 & 4.22e-03 & 8.55e-03 \\
        2 & 32 & 1185 & 7.76e-05 & 1.90e-04 & 3.05e-04 & 1.44e-03 & 2.25e-03 & 3.02e-03 \\
        2 & 64 & 4417 & 6.32e-05 & 1.71e-04 & 2.54e-04 & 8.23e-04 & 1.87e-03 & 2.69e-03 \\
        \hline
        \hline
        3 & 8 & 177 & 1.73e-04 & 3.56e-04 & 5.50e-04 & 1.22e-03 & 2.46e-03 & 6.04e-03 \\
        3 & 16 & 609 & 9.30e-05 & 2.41e-04 & 4.30e-04 & 5.94e-04 & 1.90e-03 & 3.77e-03 \\
        3 & 32 & 2241 & 6.55e-05 & 1.59e-04 & 3.71e-04 & 4.45e-04 & 1.25e-03 & 2.67e-03 \\
        3 & 64 & 8577 & 3.51e-05 & 1.74e-04 & 5.75e-04 & 2.80e-04 & 1.53e-03 & 4.29e-03 \\
        \hline
        \hline
        4 & 8 & 249 & 1.36e-04 & 2.93e-04 & 6.90e-04 & 9.37e-04 & 2.22e-03 & 3.93e-03 \\
        4 & 16 & 881 & 6.70e-05 & 1.35e-04 & 2.36e-04 & 4.31e-04 & 1.41e-03 & 2.15e-03 \\
        4 & 32 & 3297 & 5.61e-05 & 1.16e-04 & 2.82e-04 & 4.92e-04 & 1.08e-03 & 2.35e-03 \\
        4 & 64 & 12737 & 5.01e-05 & 1.02e-04 & 2.53e-04 & 3.74e-04 & 8.75e-04 & 1.97e-03 \\
        \hline
        \hline
        5 & 8 & 321 & 1.70e-04 & 2.29e-04 & 3.15e-04 & 9.16e-04 & 1.62e-03 & 2.54e-03 \\
        5 & 16 & 1153 & 8.88e-05 & 1.21e-04 & 1.45e-04 & 4.78e-04 & 1.19e-03 & 1.53e-03 \\
        5 & 32 & 4353 & 4.96e-05 & 9.23e-05 & 1.41e-04 & 5.35e-04 & 7.95e-04 & 1.32e-03 \\
        5 & 64 & 16897 & 6.12e-05 & 9.28e-05 & 1.87e-04 & 4.06e-04 & 7.16e-04 & 1.21e-03 \\
        \hline
    \end{tabular}
    \caption{Errors for Burgers equation with viscosity for MLP networks.}
    \label{tab:burgers_nn}
\end{table}

\begin{table}
    \centering
    \begin{tabular}{|c|c||c|c|}
        \hline
        Degree & d.o.f. & $\epsilon_2$ & $\epsilon_\infty$ \\
        \hline
        \hline
        2 & 81 & 1.68e+00 & 5.21e+00 \\
        2 & 289 & 5.12e-01 & 2.74e+00 \\
        2 & 1089 & 8.59e-02 & 9.12e-01 \\
        2 & 4225 & 1.26e-02 & 1.60e-01 \\
        2 & 16641 & 2.78e-03 & 3.04e-02 \\
        \hline
        \hline
        3 & 100 & 7.87e-01 & 2.97e+00 \\
        3 & 324 & 2.95e-01 & 1.90e+00 \\
        3 & 1156 & 4.41e-02 & 6.50e-01 \\
        3 & 4356 & 9.44e-03 & 1.23e-01 \\
        3 & 16900 & 2.55e-03 & 3.10e-02 \\
        \hline
        \hline
        4 & 121 & 1.54e+00 & 5.07e+00 \\
        4 & 361 & 3.40e-01 & 2.13e+00 \\
        4 & 1225 & 4.19e-02 & 6.27e-01 \\
        4 & 4489 & 1.90e-03 & 5.09e-02 \\
        4 & 17161 & 7.94e-05 & 2.19e-03 \\
        \hline
        \hline
        5 & 144 & 7.84e-01 & 2.82e+00 \\
        5 & 400 & 3.25e-01 & 2.09e+00 \\
        5 & 1296 & 3.87e-02 & 5.55e-01 \\
        5 & 4624 & 1.51e-03 & 4.47e-02 \\
        5 & 17424 & 3.75e-05 & 1.05e-03 \\
        \hline
        \hline
        6 & 169 & 1.51e+00 & 4.93e+00 \\
        6 & 441 & 3.14e-01 & 1.98e+00 \\
        6 & 1369 & 3.70e-02 & 5.46e-01 \\
        6 & 4761 & 1.24e-03 & 3.53e-02 \\
        6 & 17689 & 1.04e-05 & 4.35e-04 \\
        \hline
    \end{tabular}
    \caption{Errors for Burgers equation with viscosity for B-Splines.}
    \label{tab:burgers_sp}
\end{table}

\clearpage

\subsection{Incompressible Navier-Stokes equations}

We consider then the two-dimensional steady incompressible Navier-Stokes equations:
\begin{align}
    & u \, \partial_x u + v \, \partial_y u -\nu (\partial_{xx} u + \partial_{yy} u ) + p/\rho = 0 ,\\
    & u \, \partial_x v + v \, \partial_y v -\nu (\partial_{xx} v + \partial_{yy} v ) + p/\rho = 0 ,\\
    & \partial_x u +  \partial_y v = 0 \qquad \forall (x,y) \in ]-1,1[^2 ,
\end{align}
where $(u,v)$ are the two components of the velocity and $p$ the pressure. The density $\rho$ and the viscosity $\nu$ are supposed constant.
The main point of interest in this problem is to see how a system of equations is taken into account. When using neural networks, all unknown fields are modeled by a unique network with three outputs, which introduces a coupling between the fields in the representation. On the contrary, when using B-Splines, each unknown is represented by its own model in a fully uncoupled way, yielding an increase of the number of degrees of freedom by a factor three.

\begin{figure}[!ht]
    \centering
    \includegraphics[width=0.5\textwidth]{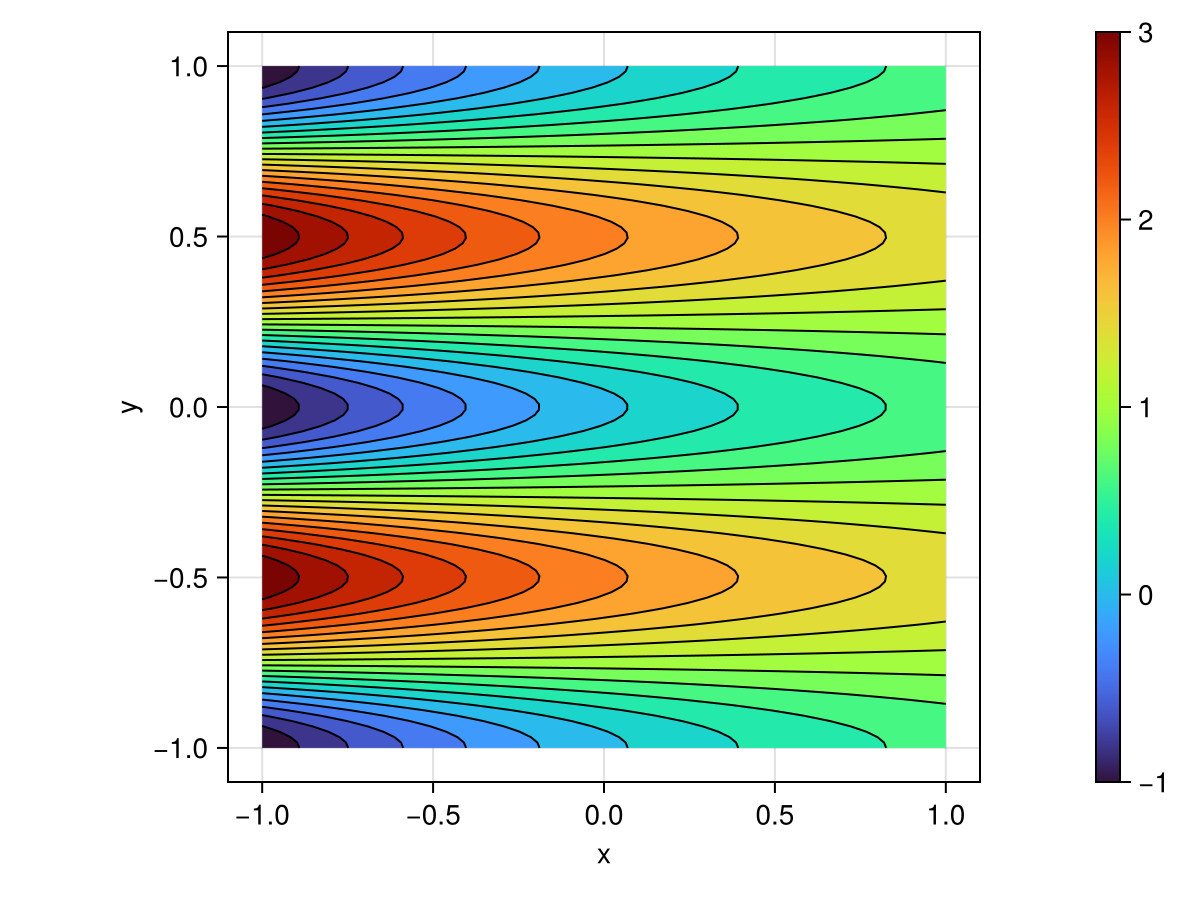}\\
    \includegraphics[width=0.5\textwidth]{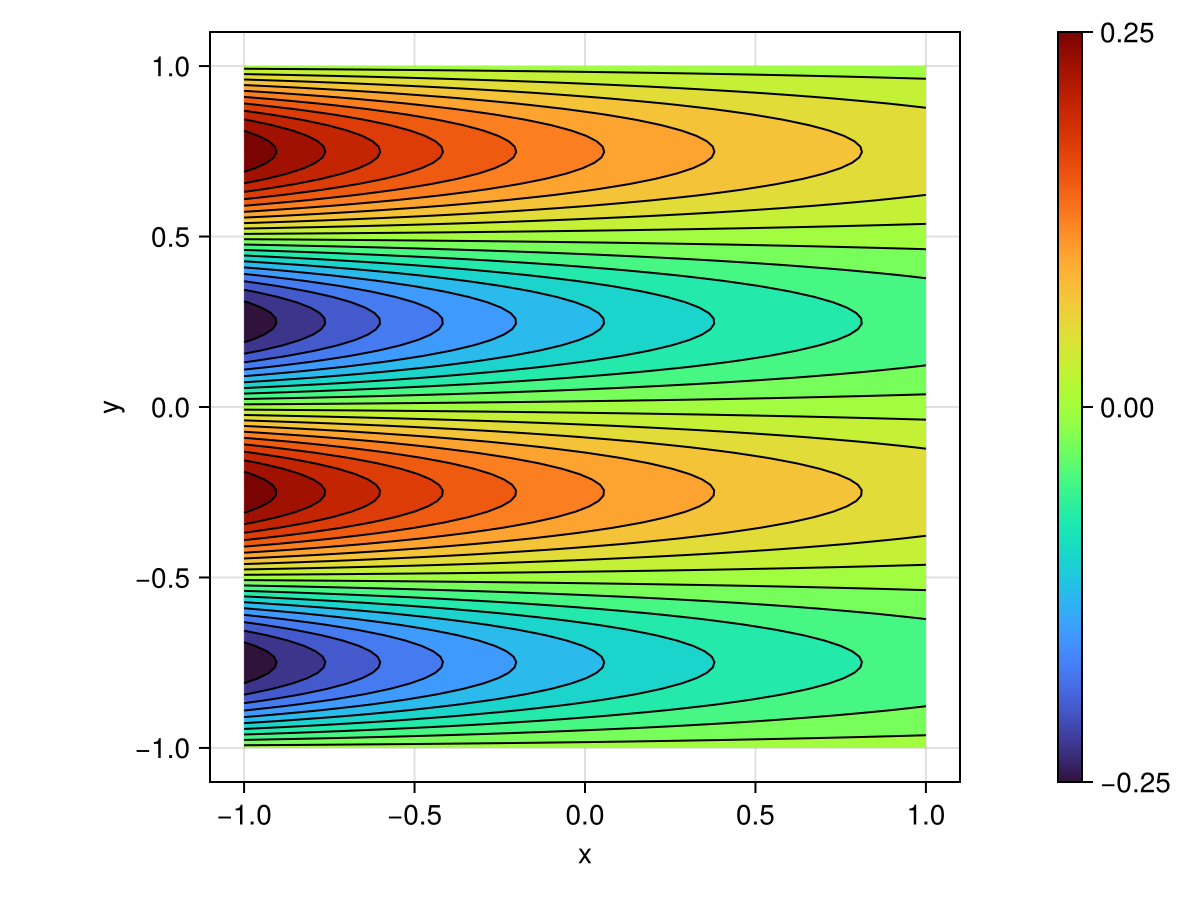}\\
    \includegraphics[width=0.5\textwidth]{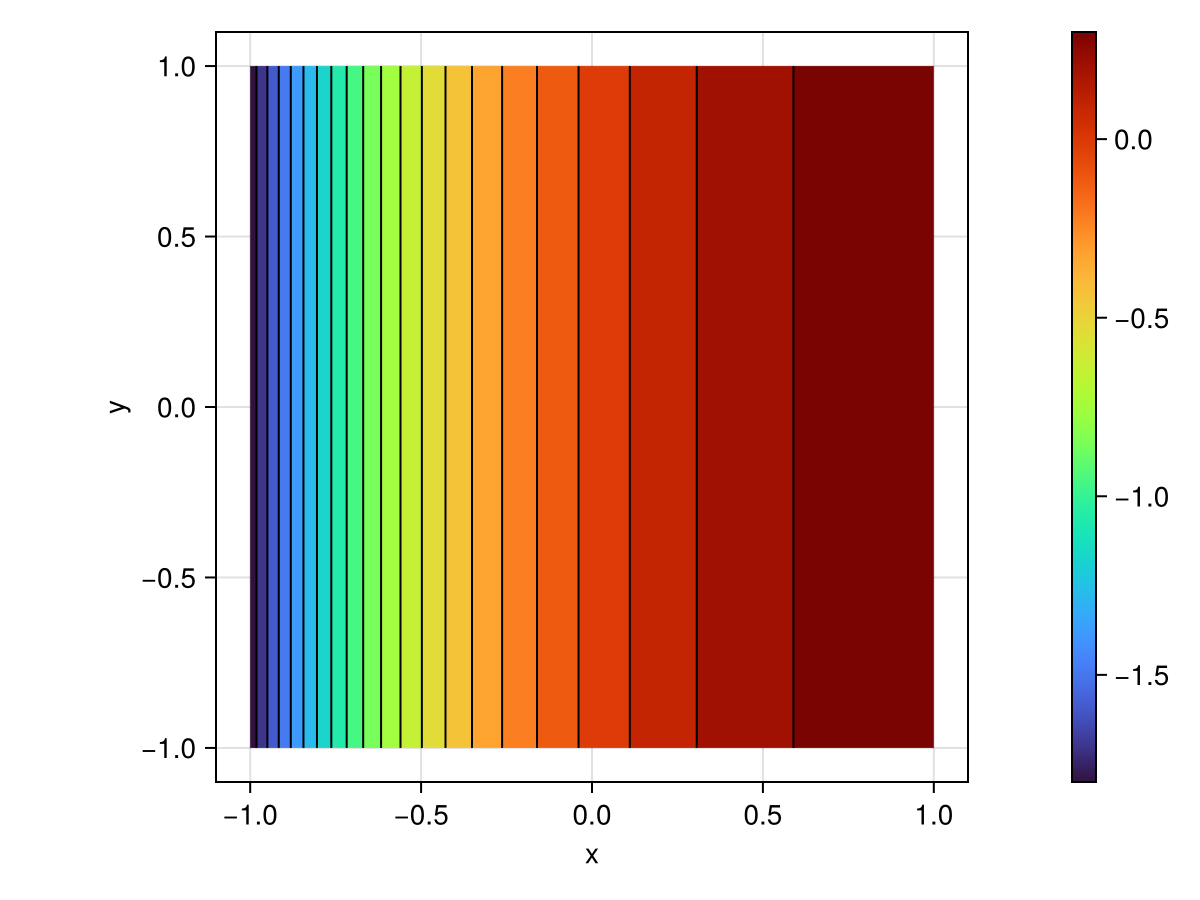}
     \caption{Analytical solution for Kovaznay problem, for $u$ (top), $v$ (middle) and $p$ (bottom).}
    \label{fig:kovaznay_exact}
\end{figure}

The closed-form solution found by Kovaznay~\cite{Kovasznay_1948} is targeted:
\begin{align}
    & u(x,y) = 1 - e^{\lambda x} \cos(2\pi y) ,\\
    & v(x,y) = \frac{\lambda}{2\pi} e^{\lambda x} \sin(2\pi y) ,\\
    & p(x,y) = \frac{1}{2}( 1 - e^{2\lambda x}) ,
\end{align}
with $\lambda = \frac{1}{2} (R - \sqrt{R^2+16\pi^2})$ and $R = UL/\nu$. This solution can be interpreted as the flow behind a two-dimensional grid. In the configuration studied, $U=L=1$, $\nu=0.02$ and $\rho=1$. The computational domain is the square $[-1,1]\times[-1,1]$. As previously, all the boundary conditions are set using the analytical solution, by projection for B-Splines and penalization for neural networks. B-Spline parameters are initialized with zero value, except at boundaries. B-Spline representations are tested using 8, 16, 32 and 64 parameters in each direction and each unknown field. The exact solution is shown on Fig.~(\ref{fig:kovaznay_exact}).

\begin{figure}[!ht]
  \begin{subfigure}[b]{0.49\linewidth}
    \centering
    \includegraphics[width=\textwidth]{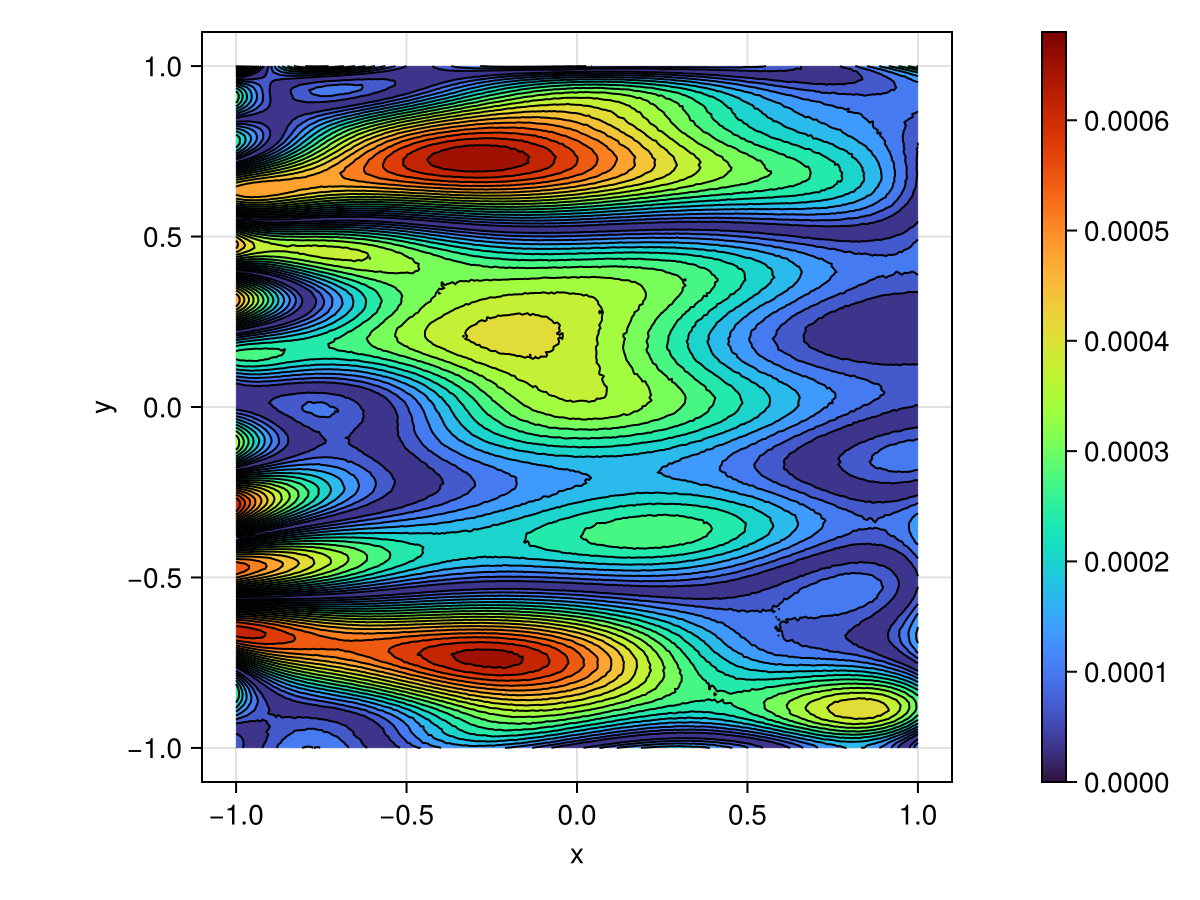}
     \caption{Neural network}
    \label{fig:kozvanay_error_nn}
  \end{subfigure}
  \begin{subfigure}[b]{0.49\linewidth}
    \centering
    \includegraphics[width=\textwidth]{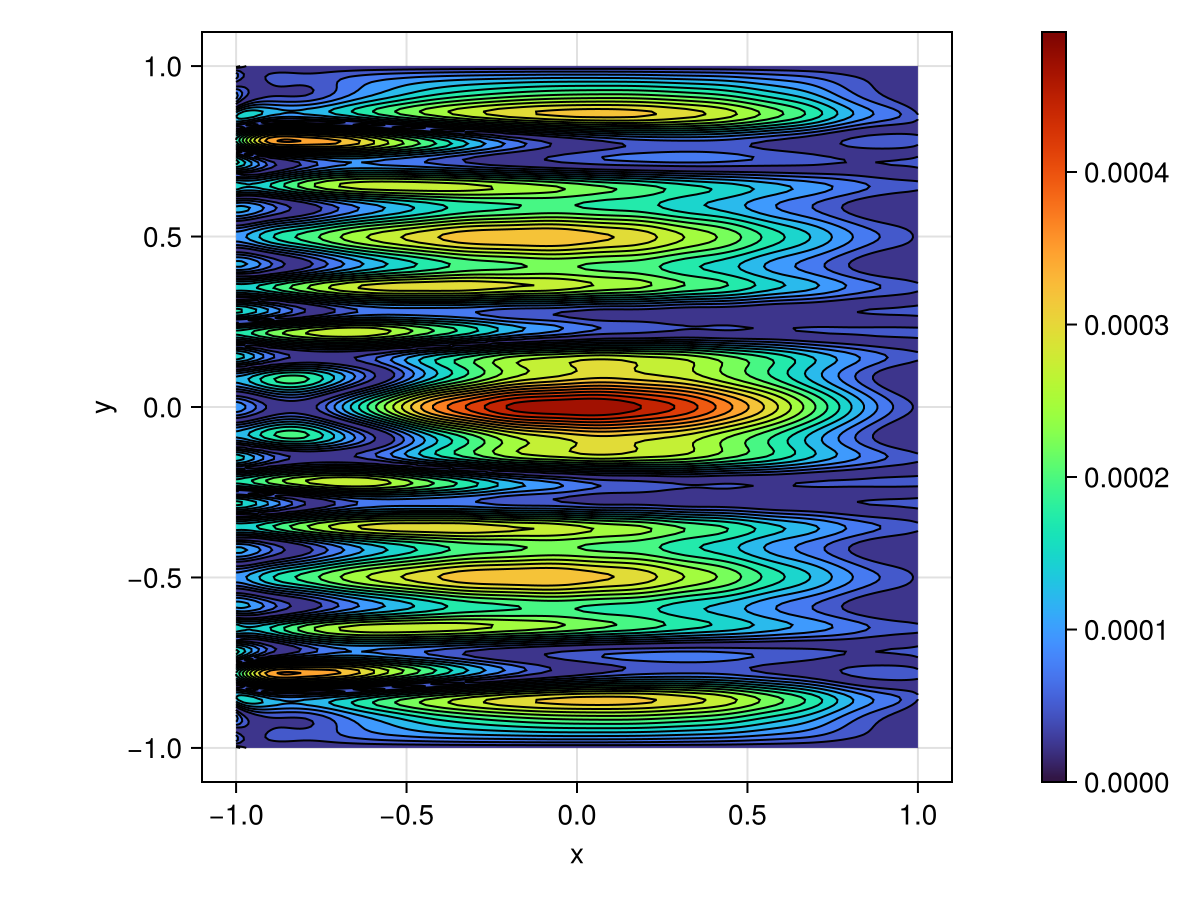}
    \caption{B-Spline basis}
    \label{fig:kozvanay_error_sp}
  \end{subfigure}
  \caption{Error fields for $u$ for Navier-Stokes equations (about 1200 parameters).}
  \label{fig:kozvanay_error}
\end{figure}

\begin{figure}[!ht]
  \begin{subfigure}[b]{0.49\linewidth}
    \centering
    \includegraphics[width=\textwidth]{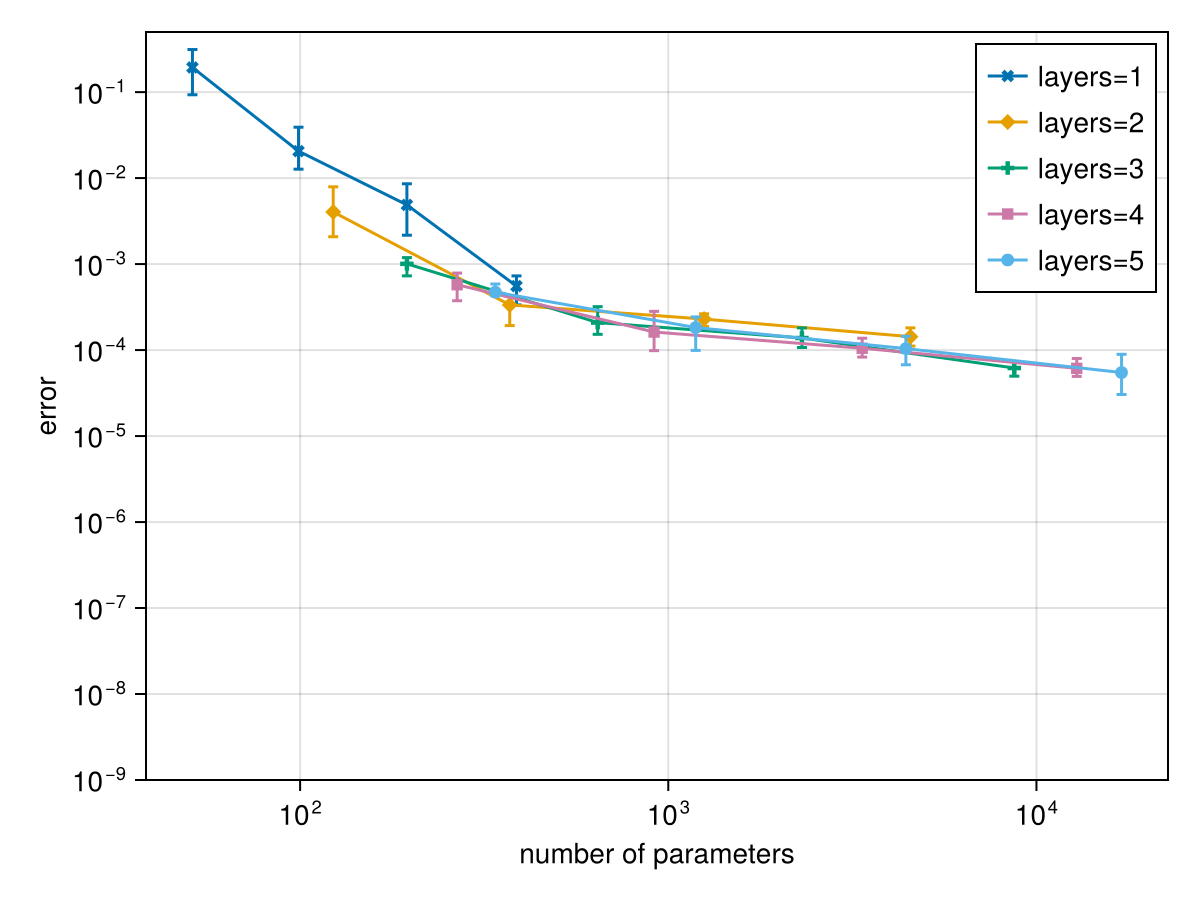}
     \caption{Neural network}
    \label{fig:kovaznay_nn_l2}
  \end{subfigure}
  \begin{subfigure}[b]{0.49\linewidth}
    \centering
    \includegraphics[width=\textwidth]{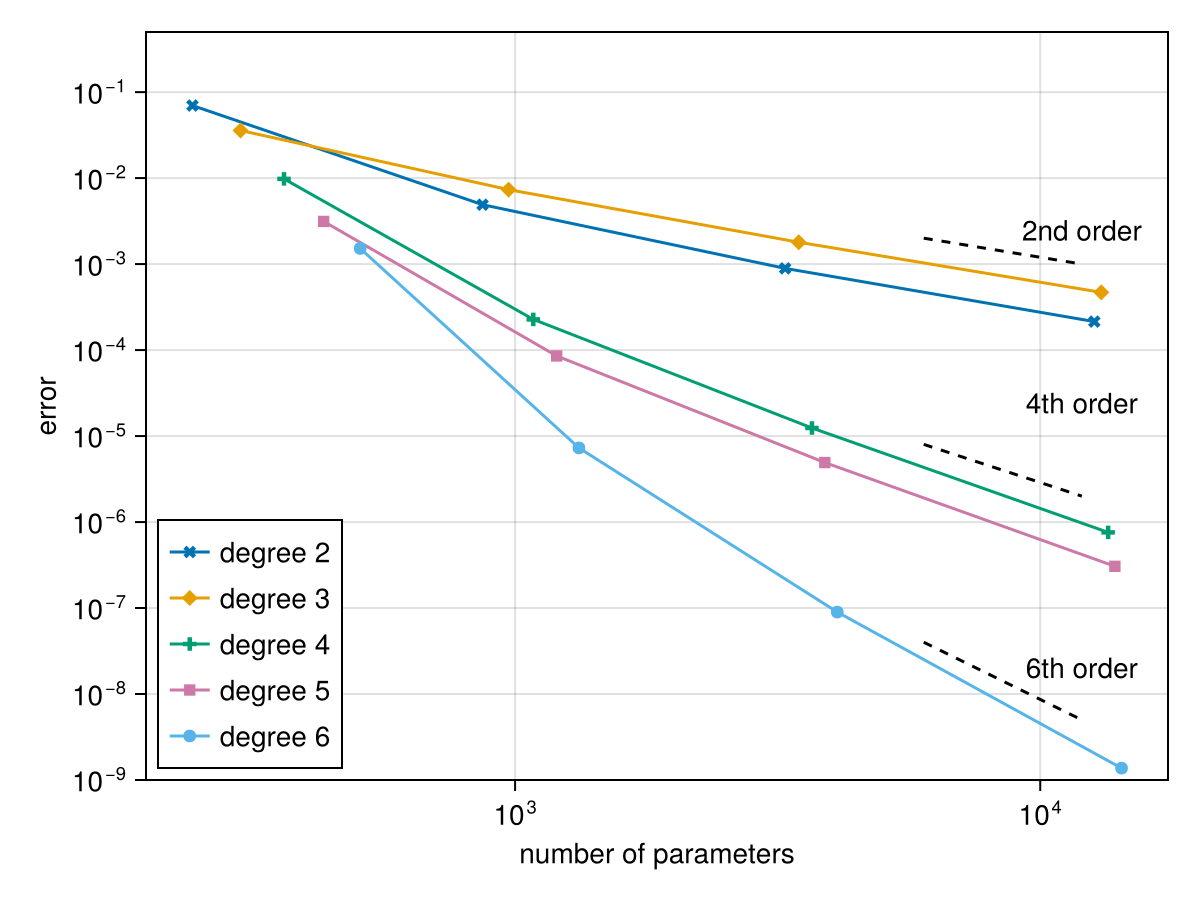}
    \caption{Spline basis}
    \label{fig:kovaznay_sp_l2}
  \end{subfigure}
  \caption{$\epsilon_2$ norm of error w.r.t. the number of parameters for Navier-Stokes equations.}
  \label{fig:kovaznay_l2}
\end{figure}

The variation in error as a function of the number of parameters is depicted in Fig.~(\ref{fig:kovaznay_l2}), for the field $u$, while error values are provided in Tabs.~(\ref{tab:ns_nn}-\ref{tab:ns_sp}). For B-Splines, the total number of parameters is reported, accounting for the three unknown fields. We find trends similar to those observed for the elliptic problem, except for the single-layer networks  which exhibit higher error levels. Regarding B-Splines, the theoretical asymptotic convergence rates are perfectly recovered.
The error maps in Fig.~(\ref{fig:kozvanay_error}) show a perfectly symmetric error field for B-Splines, whereas networks cannot maintain the symetry of the solution due to the randomness of the sampling and initialization.

\begin{table}
    \centering
    \begin{tabular}{|c|c|c||c|c|c|c|c|c| }
        \hline
        Layers & Neurons & d.o.f. & \multicolumn{3}{|c|}{$\epsilon_2$} & \multicolumn{3}{|c|}{$\epsilon_\infty$} \\
        \hline
         & & & min & average & max & min & average & max \\
        \hline
        \hline
        1 & 8 & 51 & 9.31e-02 & 1.94e-01 & 3.13e-01 & 3.45e-01 & 5.12e-01 & 1.03e+00 \\
        1 & 16 & 99 & 1.27e-02 & 2.06e-02 & 3.91e-02 & 3.19e-02 & 6.93e-02 & 1.22e-01 \\
        1 & 32 & 195 & 2.17e-03 & 4.88e-03 & 8.59e-03 & 1.22e-02 & 2.11e-02 & 3.75e-02 \\
        1 & 64 & 387 & 3.35e-04 & 5.54e-04 & 7.28e-04 & 2.48e-03 & 4.34e-03 & 5.99e-03 \\
        \hline
        \hline
        2 & 8 & 123 & 2.08e-03 & 4.04e-03 & 7.91e-03 & 5.85e-03 & 1.41e-02 & 2.85e-02 \\
        2 & 16 & 371 & 1.93e-04 & 3.35e-04 & 4.22e-04 & 8.30e-04 & 1.35e-03 & 1.68e-03 \\
        2 & 32 & 1251 & 1.88e-04 & 2.29e-04 & 2.64e-04 & 6.78e-04 & 9.64e-04 & 1.35e-03 \\
        2 & 64 & 4547 & 1.11e-04 & 1.43e-04 & 1.81e-04 & 6.14e-04 & 7.23e-04 & 9.03e-04 \\
        \hline
        \hline
        3 & 8 & 195 & 7.29e-04 & 1.01e-03 & 1.19e-03 & 3.02e-03 & 3.75e-03 & 4.63e-03 \\
        3 & 16 & 643 & 1.52e-04 & 2.09e-04 & 3.20e-04 & 7.07e-04 & 9.15e-04 & 1.16e-03 \\
        3 & 32 & 2307 & 1.07e-04 & 1.38e-04 & 1.81e-04 & 6.16e-04 & 7.05e-04 & 8.82e-04 \\
        3 & 64 & 8707 & 4.97e-05 & 6.18e-05 & 7.35e-05 & 2.23e-04 & 4.12e-04 & 6.08e-04 \\
        \hline
        \hline
        4 & 8 & 267 & 3.75e-04 & 5.75e-04 & 7.86e-04 & 1.49e-03 & 2.33e-03 & 3.44e-03 \\
        4 & 16 & 915 & 9.84e-05 & 1.62e-04 & 2.82e-04 & 5.10e-04 & 7.46e-04 & 1.22e-03 \\
        4 & 32 & 3363 & 8.29e-05 & 1.05e-04 & 1.37e-04 & 4.36e-04 & 4.78e-04 & 5.37e-04 \\
        4 & 64 & 12867 & 4.94e-05 & 6.13e-05 & 7.98e-05 & 3.33e-04 & 3.82e-04 & 4.79e-04 \\
        \hline
        \hline
        5 & 8 & 339 & 4.27e-04 & 4.71e-04 & 5.87e-04 & 1.17e-03 & 1.55e-03 & 2.24e-03 \\
        5 & 16 & 1187 & 9.91e-05 & 1.83e-04 & 2.43e-04 & 5.86e-04 & 8.66e-04 & 1.10e-03 \\
        5 & 32 & 4419 & 6.74e-05 & 1.04e-04 & 1.45e-04 & 4.21e-04 & 5.41e-04 & 7.92e-04 \\
        5 & 64 & 17027 & 3.04e-05 & 5.47e-05 & 8.93e-05 & 1.37e-04 & 3.78e-04 & 7.68e-04 \\
        \hline
    \end{tabular}
    \caption{Errors for Navier-Stokes equations for MLP networks ($u$ field).}
    \label{tab:ns_nn}
\end{table}

\begin{table}
    \centering
    \begin{tabular}{|c|c||c|c|}
        \hline
        Degree & d.o.f. & $\epsilon_2$ & $\epsilon_\infty$ \\
        \hline
        \hline
        2 & 243 & 6.98e-02 & 1.91e-01 \\
        2 & 867 & 4.90e-03 & 1.39e-02 \\
        2 & 3267 & 8.91e-04 & 2.35e-03 \\
        2 & 12675 & 2.14e-04 & 5.66e-04 \\
        \hline
        \hline
        3 & 300 & 3.58e-02 & 9.19e-02 \\
        3 & 972 & 7.33e-03 & 2.18e-02 \\
        3 & 3468 & 1.79e-03 & 5.34e-03 \\
        3 & 13068 & 4.69e-04 & 1.33e-03 \\
        \hline
        \hline
        4 & 363 & 9.79e-03 & 3.32e-02 \\
        4 & 1083 & 2.27e-04 & 5.93e-04 \\
        4 & 3675 & 1.24e-05 & 3.50e-05 \\
        4 & 13467 & 7.58e-07 & 2.14e-06 \\
        \hline
        \hline
        5 & 432 & 3.13e-03 & 1.16e-02 \\
        5 & 1200 & 8.54e-05 & 2.79e-04 \\
        5 & 3888 & 4.92e-06 & 1.48e-05 \\
        5 & 13872 & 3.06e-07 & 8.74e-07 \\
        \hline
        \hline
        6 & 507 & 1.52e-03 & 3.96e-03 \\
        6 & 1323 & 7.27e-06 & 1.74e-05 \\
        6 & 4107 & 8.97e-08 & 2.53e-07 \\
        6 & 14283 & 1.37e-09 & 5.31e-09 \\
        \hline
    \end{tabular}
    \caption{Errors for Navier-Stokes equations for B-Splines ($u$ field).}
    \label{tab:ns_sp}
\end{table}

\clearpage

\subsection{Compressible Euler equations}

\begin{figure}[!ht]
  \begin{subfigure}[b]{0.49\linewidth}
    \centering
    \includegraphics[width=\textwidth]{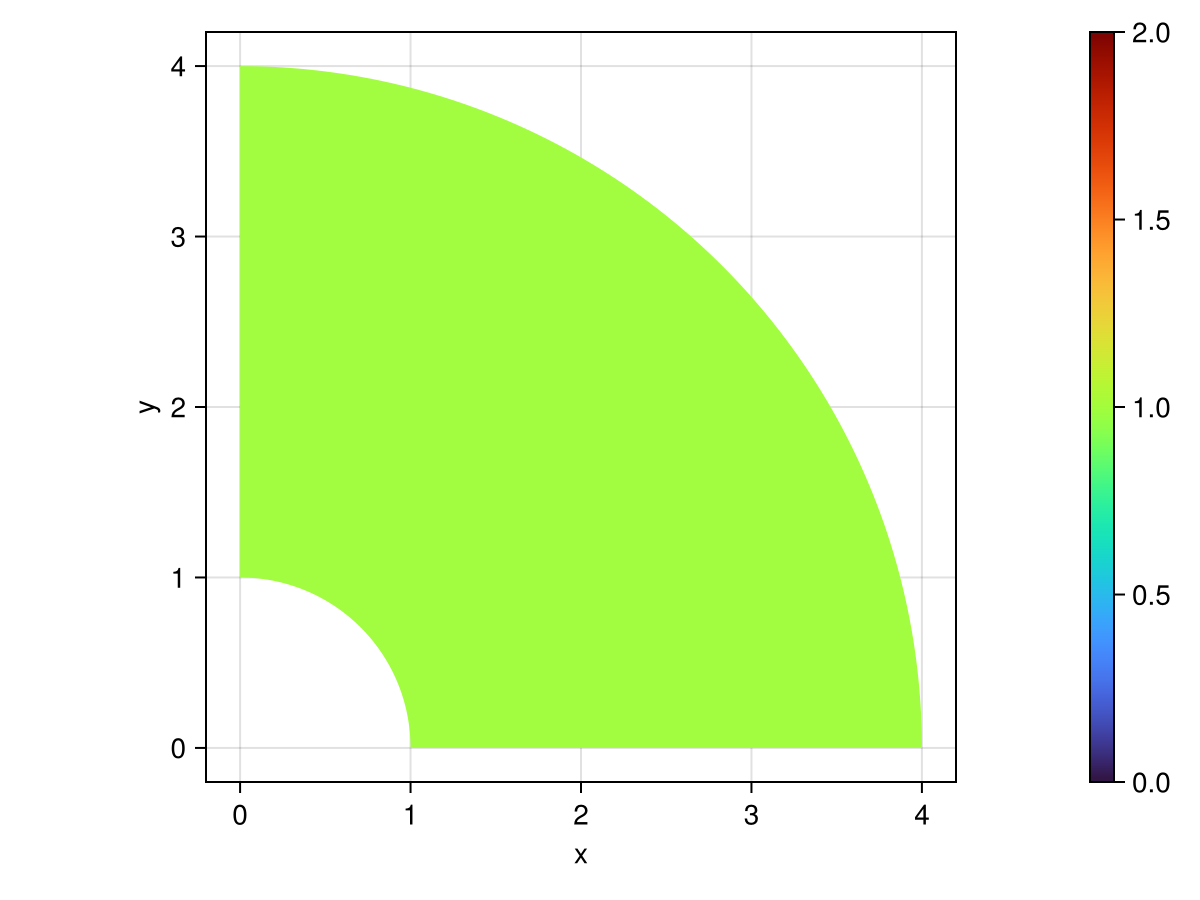}
     \caption{density}
    \label{fig:rho_exact}
  \end{subfigure}
  \begin{subfigure}[b]{0.49\linewidth}
    \centering
    \includegraphics[width=\textwidth]{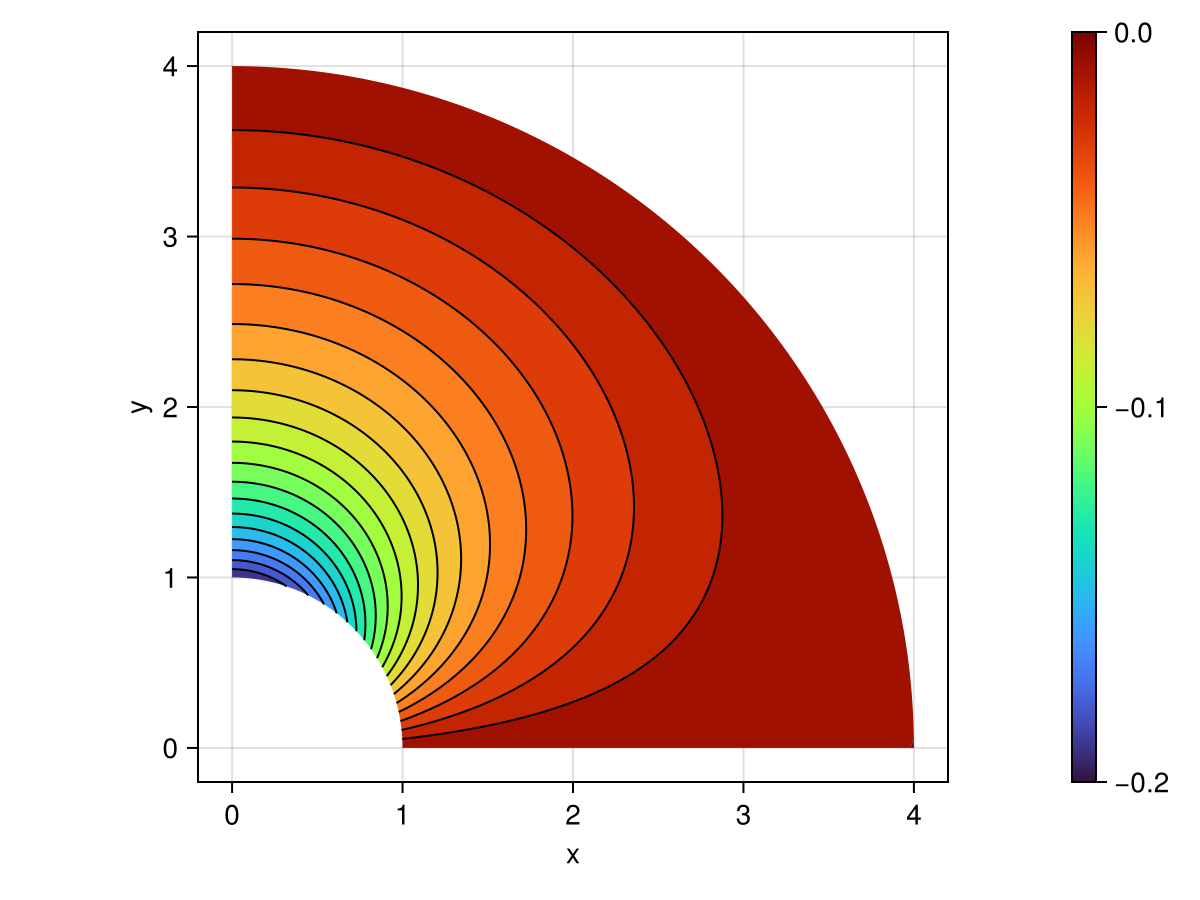}
    \caption{velocity along $x$}
    \label{fig:u_exact}
  \end{subfigure}\\
  \begin{subfigure}[b]{0.49\linewidth}
    \centering
    \includegraphics[width=\textwidth]{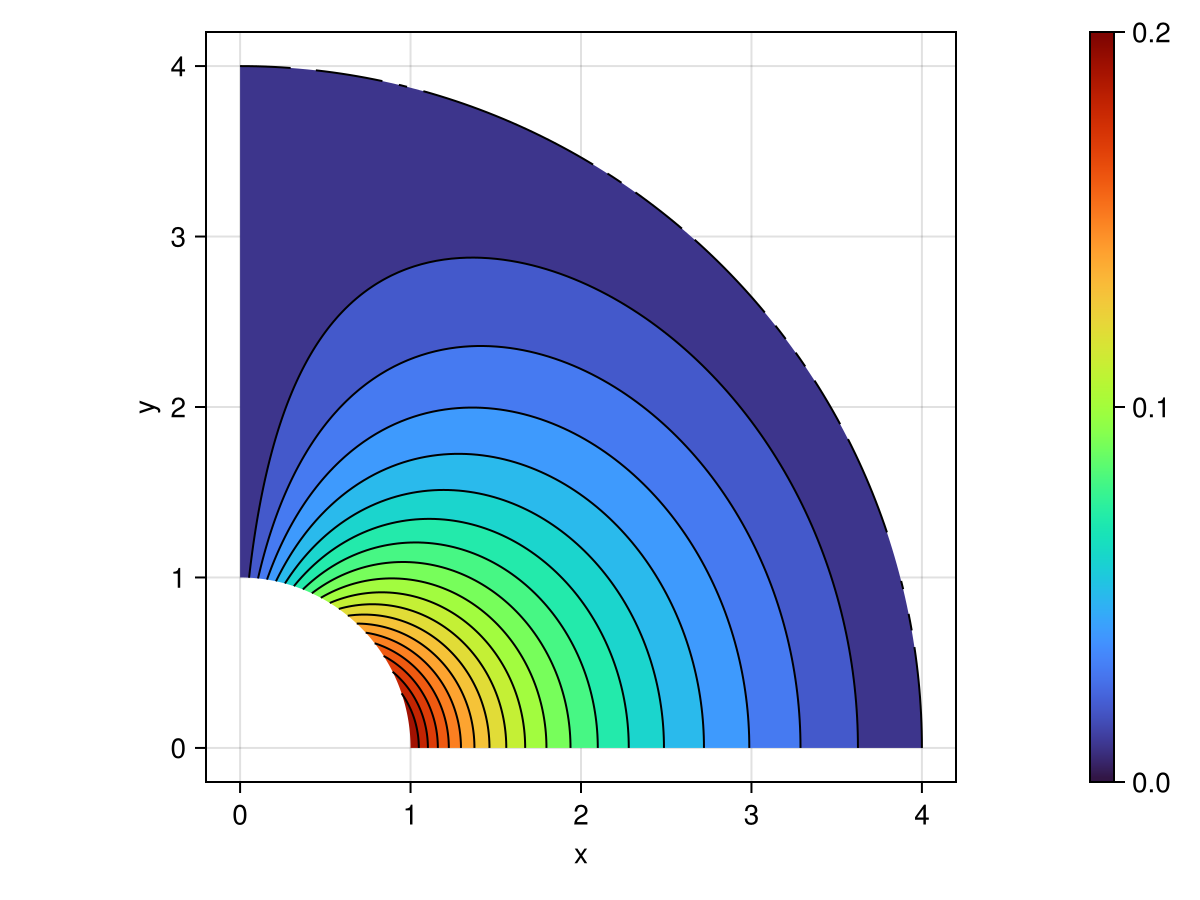}
     \caption{velocity along $y$}
    \label{fig:v_exact}
  \end{subfigure}
  \begin{subfigure}[b]{0.49\linewidth}
    \centering
    \includegraphics[width=\textwidth]{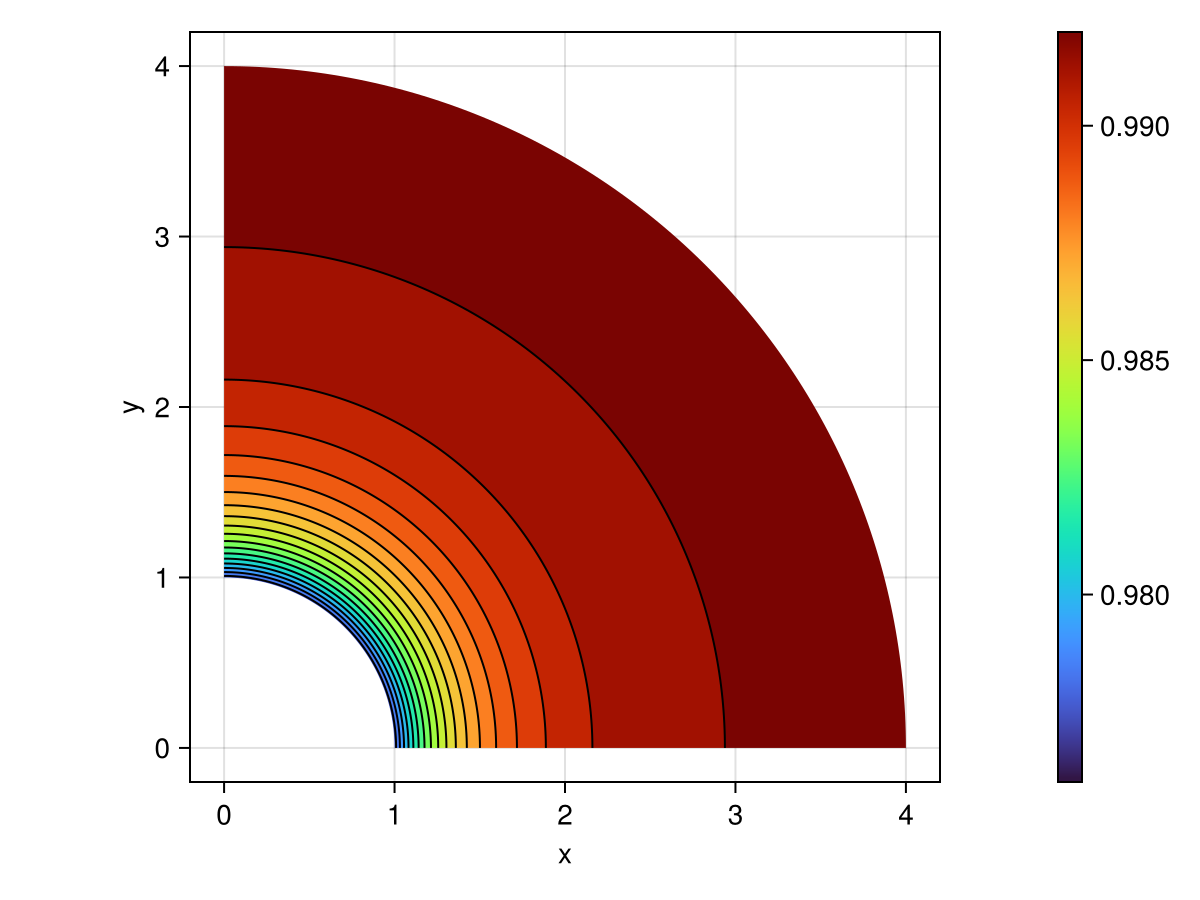}
    \caption{pressure}
    \label{fig:p_exact}
  \end{subfigure}
  \caption{Analytical solution for compressible Euler equations.}
  \label{fig:euler_exact}
\end{figure}

For the last case, the non-linear system of steady compressible Euler equations is considered:
\begin{align}
    & \partial_x (\rho u) + \partial_y (\rho v) = 0 \\
    & \partial_x (\rho u^2+p) + \partial_y (\rho u v) = 0 \\
    & \partial_x (\rho v u) + \partial_y (\rho v^2+p) = 0 \\
    & \partial_x (u ( \rho E + p )) + \partial_y (v ( \rho E + p )) = 0
\end{align}
where $\rho$ is the density, $(u,v)$ the velocity components, $E$ the total energy
and $p = \rho (\gamma - 1) ( E - \frac{1}{2} (u^2 + v^2) )$ the pressure. Note that there are no diffusion terms and the unsteady version of this system is purely hyperbolic. We consider here as testcase the flow between two cylinders, whose closed-form solution is depicted in Fig.~(\ref{fig:euler_exact}) and writes:
\begin{align}
    & \rho = 1 , \\
    & u = -\sin (\theta) \; u_\theta , \\
    & v = \cos (\theta) \; u_\theta , \\
    & p = 1 +  \frac{1}{75^2} \left ( \frac{r^2}{2} - 32 \, ln(r) - \frac{128}{r^2} \right ) ,
\end{align}
with $u_\theta = \frac{1}{75} (-r + \frac{16}{r})$. Contrary to the previous testcases, the computational domain is not a square, which impacts the learning process. This does not modify the implementation of the approach based on neural networks. Nevertheless, boundary conditions are now imposed on curved boundaries that do not correspond to constant $x$ or $y$ values, which makes the training more difficult. For B-Splines, a mapping is first defined between the parametric domain and the physical one by solving a least-square fitting problem, as commonly done in isogeometric analysis~\cite{Duvigneau_18}, yielding a relationship between the parametric coordinates and the Cartesian coordinates $x(\xi_1,\xi_2)$ and $y(\xi_1,\xi_2)$. This mapping is used to compute analytically the derivatives, e.g. $\partial_x u = \partial_{\xi_1} u \, \partial_x \xi_1 + \partial_{\xi_2} u \, \partial_x \xi_2$.
Thanks to this mapping, the boundaries of the computational domain still correspond to B-Spline extremities.

The inlet and outlet boundary conditions, imposed at $x=0$ and $y=0$ boundaries, are provided by the analytical solution of the problem. On the two curved boundaries, only a slip boundary condition is prescribed for the velocity, i.e. $\vec{U} \cdot \vec{n}$. All boundary conditions are taken into account by penalization for neural networks, whereas B-Splines rely on projection methods for inlet and outlet conditions associated with penalization for the slip condition, which is not known beforehand.

\begin{figure}[!ht]
  \begin{subfigure}[b]{0.49\linewidth}
    \centering
    \includegraphics[width=\textwidth]{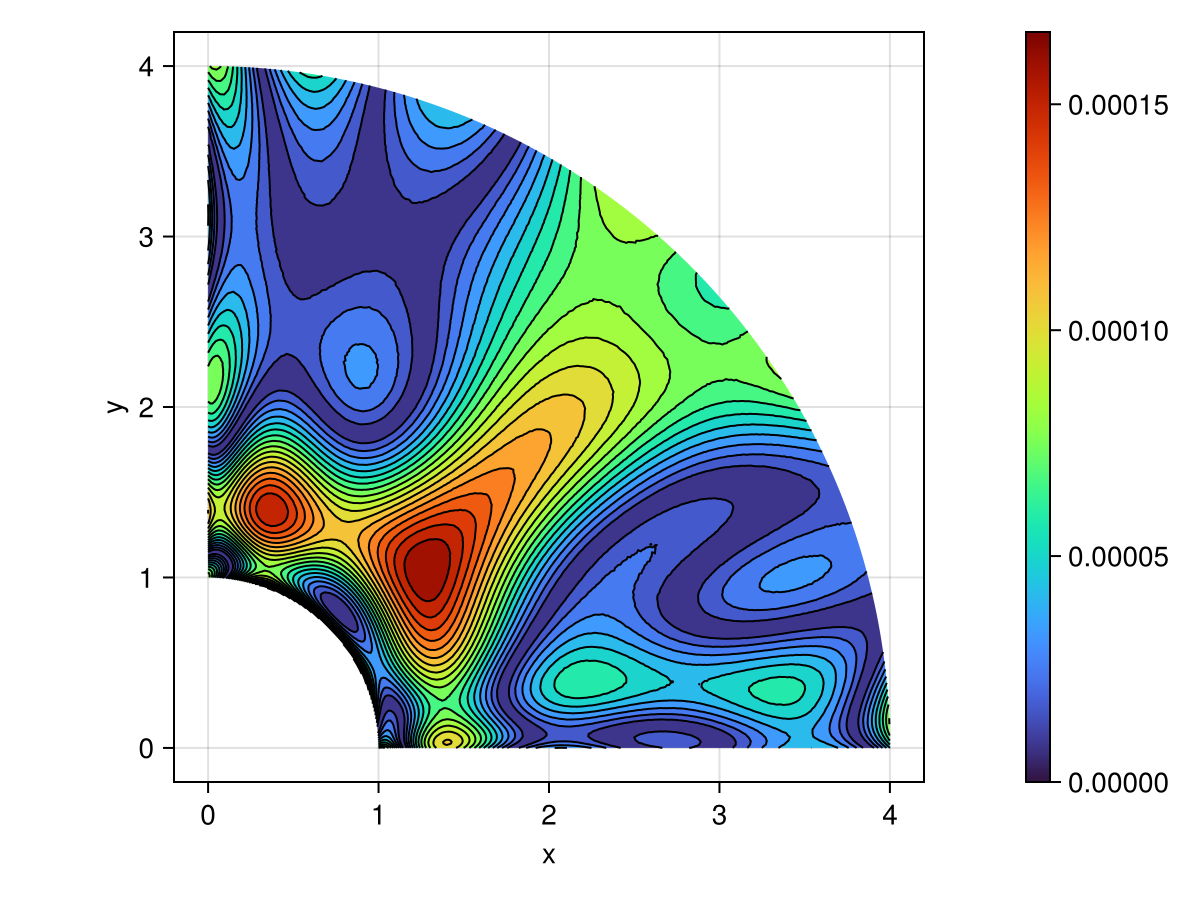}
     \caption{Neural network}
    \label{fig:euler_error_nn}
  \end{subfigure}
  \begin{subfigure}[b]{0.49\linewidth}
    \centering
    \includegraphics[width=\textwidth]{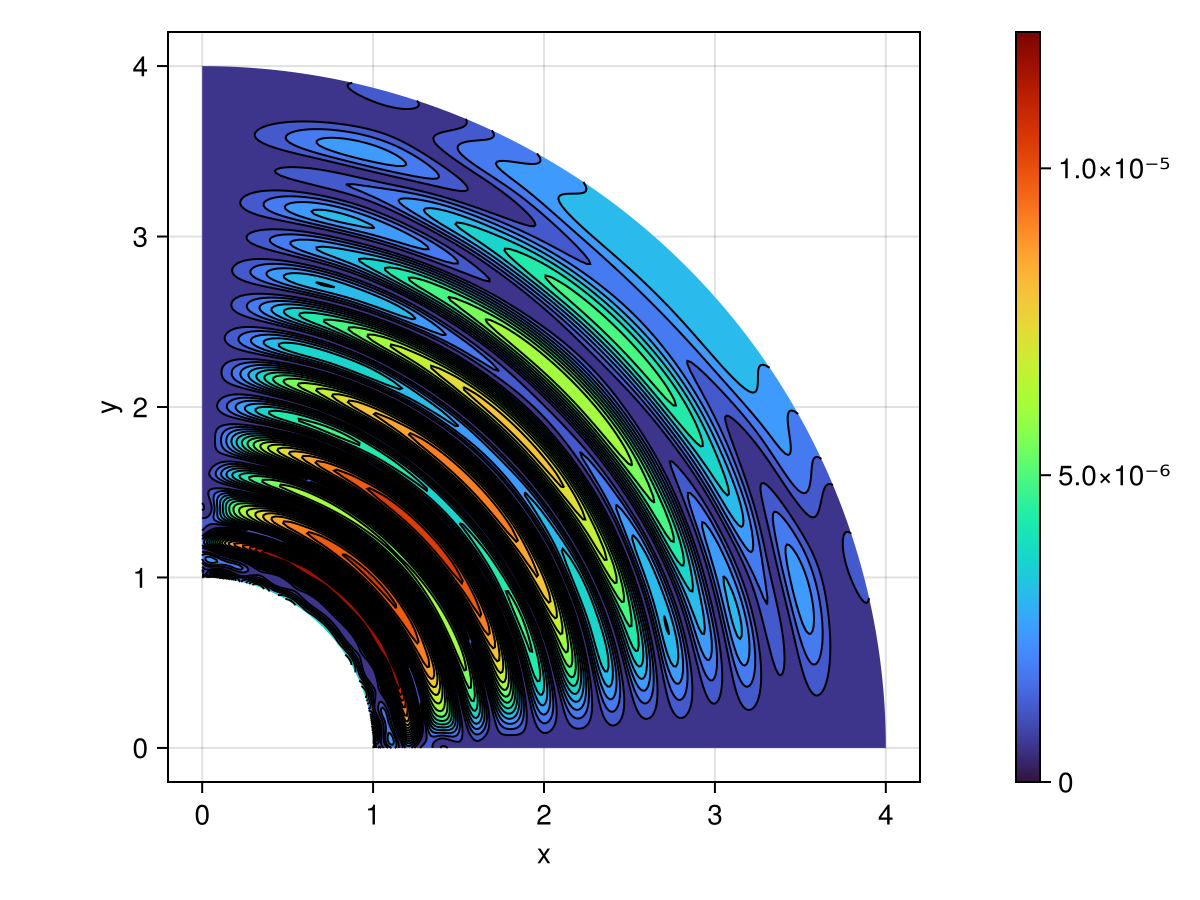}
    \caption{B-Spline basis}
    \label{fig:euler_error_sp}
  \end{subfigure}
  \caption{Error fields for $p$ for Euler equations (about 1200 parameters).}
  \label{fig:euler_error}
\end{figure}

\begin{figure}[!ht]
  \begin{subfigure}[b]{0.49\linewidth}
    \centering
    \includegraphics[width=\textwidth]{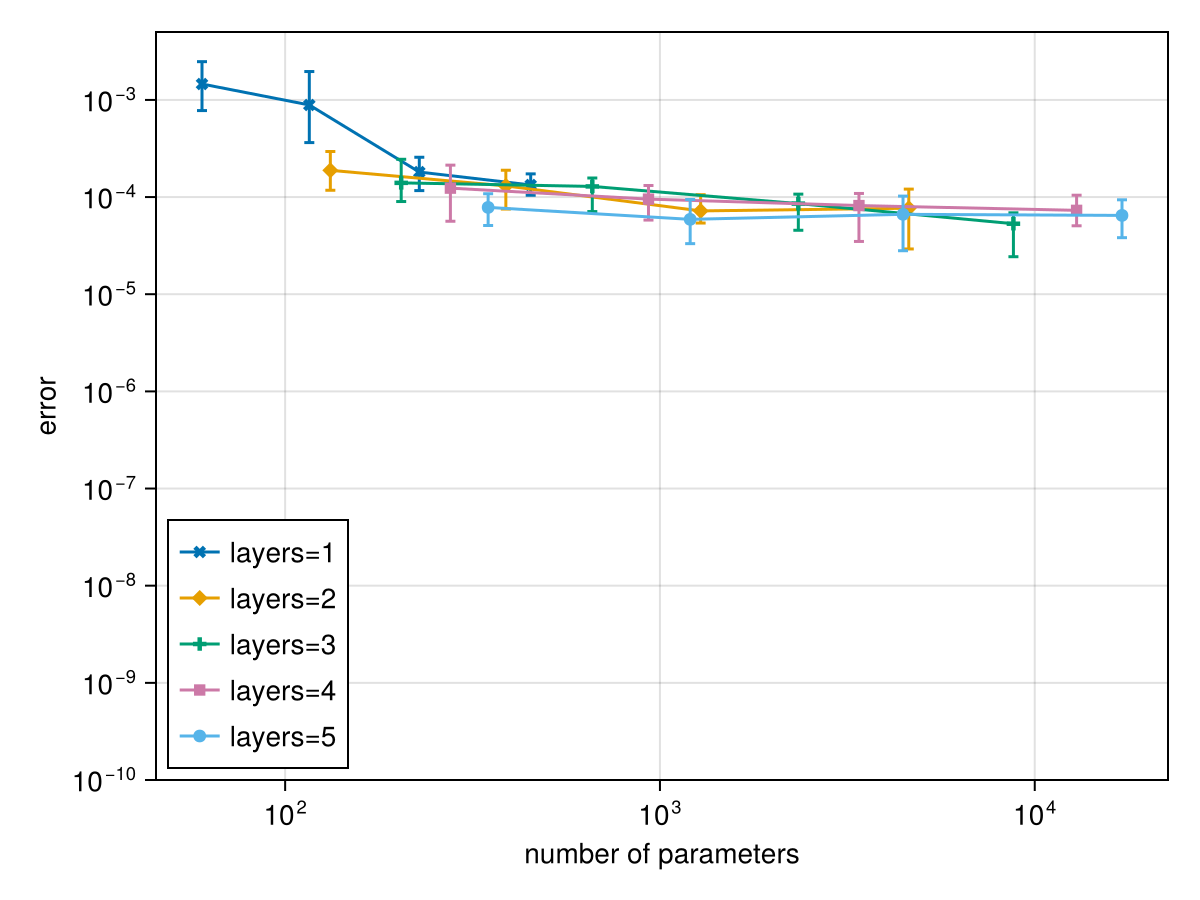}
     \caption{Neural network}
    \label{fig:euler_nn_wall_l2}
  \end{subfigure}
  \begin{subfigure}[b]{0.49\linewidth}
    \centering
    \includegraphics[width=\textwidth]{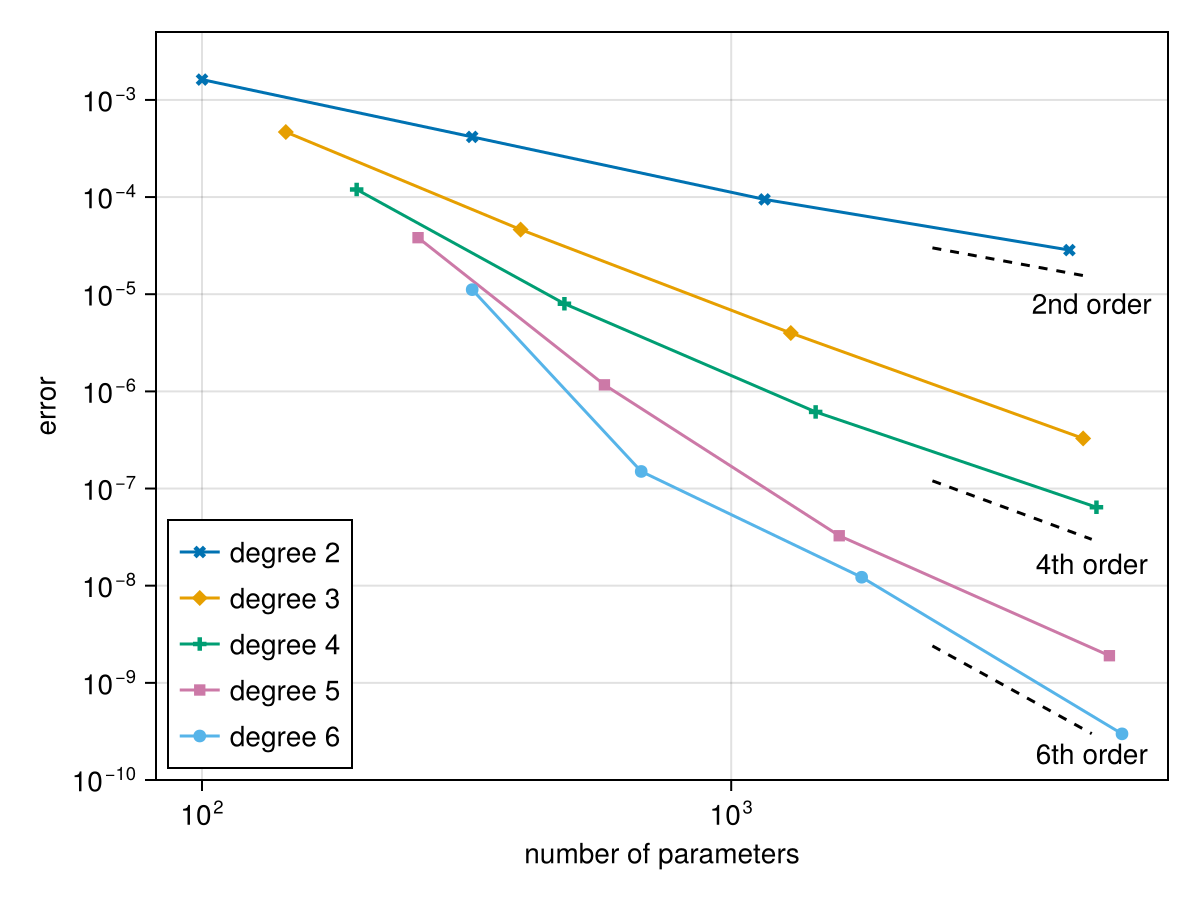}
    \caption{Spline basis}
    \label{fig:euler_sp_wall_l2}
  \end{subfigure}
  \caption{$\epsilon_2$ norm of error w.r.t. the number of parameters for $p$ field for Euler equations.}
  \label{fig:euler_wall_l2}
\end{figure}

The variation in error as a function of the number of parameters is depicted in Fig.~(\ref{fig:euler_wall_l2}) in log scale and values are provided in Tabs.~(\ref{tab:euler_nn}-\ref{tab:euler_sp}), for the pressure field $p$. As previously, one can observe a significant difference between the results obtained by the two representations. While the B-Splines exhibit an error tending to zero as the number of parameters increases and achieve a reduction of the error in agreement with the theoretical asymptotic convergence rates, the neural networks fail in reducing the error below $5 \, 10^{-5}$, whatever the architecture chosen. As shown in Fig.~(\ref{fig:euler_error}), for a number of parameters close to one thousand and two hundreds, B-Splines allow to reach an error level below $10^{-5}$, whereas the error for neural networks remains around $10^{-4}$.
Again, the symmetry is lost by solving the problem with neural networks, contrary to B-Splines.

\begin{table}
    \centering
    \begin{tabular}{|c|c|c||c|c|c|c|c|c| }
        \hline
        Layers & Neurons & d.o.f. & \multicolumn{3}{|c|}{$\epsilon_2$} & \multicolumn{3}{|c|}{$\epsilon_\infty$} \\
        \hline
         & & & min & average & max & min & average & max \\
        \hline
        \hline
        1 & 8 & 60 & 7.77e-04 & 1.46e-03 & 2.48e-03 & 3.01e-03 & 4.59e-03 & 9.37e-03 \\
        1 & 16 & 116 & 3.64e-04 & 8.87e-04 & 1.96e-03 & 1.28e-03 & 2.47e-03 & 4.77e-03 \\
        1 & 32 & 228 & 1.17e-04 & 1.81e-04 & 2.57e-04 & 4.02e-04 & 5.29e-04 & 6.64e-04 \\
        1 & 64 & 452 & 1.04e-04 & 1.33e-04 & 1.73e-04 & 3.32e-04 & 3.83e-04 & 4.35e-04 \\
        \hline
        \hline
        2 & 8 & 132 & 1.17e-04 & 1.88e-04 & 2.94e-04 & 3.38e-04 & 6.42e-04 & 1.02e-03 \\
        2 & 16 & 388 & 7.53e-05 & 1.30e-04 & 1.89e-04 & 1.87e-04 & 3.76e-04 & 6.93e-04 \\
        2 & 32 & 1284 & 5.40e-05 & 7.20e-05 & 1.06e-04 & 1.39e-04 & 1.90e-04 & 2.91e-04 \\
        2 & 64 & 4612 & 2.92e-05 & 7.68e-05 & 1.21e-04 & 1.68e-04 & 2.16e-04 & 2.72e-04 \\
        \hline
        \hline
        3 & 8 & 204 & 8.99e-05 & 1.40e-04 & 2.45e-04 & 2.85e-04 & 4.24e-04 & 6.03e-04 \\
        3 & 16 & 660 & 7.09e-05 & 1.29e-04 & 1.57e-04 & 2.22e-04 & 3.23e-04 & 3.73e-04 \\
        3 & 32 & 2340 & 4.55e-05 & 8.56e-05 & 1.07e-04 & 1.27e-04 & 2.03e-04 & 2.71e-04 \\
        3 & 64 & 8772 & 2.43e-05 & 5.32e-05 & 6.88e-05 & 1.28e-04 & 1.60e-04 & 1.80e-04 \\
        \hline
        \hline
        4 & 8 & 276 & 5.64e-05 & 1.24e-04 & 2.13e-04 & 2.22e-04 & 3.04e-04 & 4.54e-04 \\
        4 & 16 & 932 & 5.80e-05 & 9.52e-05 & 1.32e-04 & 1.56e-04 & 2.30e-04 & 2.98e-04 \\
        4 & 32 & 3396 & 3.49e-05 & 8.18e-05 & 1.09e-04 & 1.01e-04 & 1.89e-04 & 2.40e-04 \\
        4 & 64 & 12932 & 5.06e-05 & 7.30e-05 & 1.04e-04 & 1.20e-04 & 1.68e-04 & 2.50e-04 \\
        \hline
        \hline
        5 & 8 & 348 & 5.10e-05 & 7.82e-05 & 1.08e-04 & 1.71e-04 & 2.30e-04 & 3.51e-04 \\
        5 & 16 & 1204 & 3.32e-05 & 5.90e-05 & 9.44e-05 & 8.64e-05 & 1.59e-04 & 2.10e-04 \\
        5 & 32 & 4452 & 2.80e-05 & 6.64e-05 & 1.02e-04 & 1.13e-04 & 1.68e-04 & 2.57e-04 \\
        5 & 64 & 17092 & 3.82e-05 & 6.47e-05 & 9.37e-05 & 9.72e-05 & 1.56e-04 & 2.06e-04 \\
        \hline
    \end{tabular}
    \caption{Errors for Euler equations for MLP networks ($p$ field).}
    \label{tab:euler_nn}
\end{table}

\begin{table}
    \centering
    \begin{tabular}{|c|c||c|c|}
        \hline
        Degree & d.o.f. & $\epsilon_2$ & $\epsilon_\infty$ \\
        \hline
        \hline
        2 & 100 & 1.62e-03 & 4.78e-03 \\
        2 & 364 & 4.16e-04 & 1.24e-03 \\
        2 & 1156 & 9.46e-05 & 2.23e-04 \\
        2 & 4356 & 2.85e-05 & 1.33e-04 \\
        \hline
        \hline
        3 & 144 & 4.67e-04 & 1.38e-03 \\
        3 & 400 & 4.63e-05 & 1.26e-04 \\
        3 & 1296 & 3.99e-06 & 1.21e-05 \\
        3 & 4624 & 3.28e-07 & 9.28e-07 \\
        \hline
        \hline
        4 & 196 & 1.20e-04 & 2.69e-04 \\
        4 & 484 & 7.98e-06 & 1.87e-05 \\
        4 & 1444 & 6.14e-07 & 1.46e-06 \\
        4 & 4900 & 6.42e-08 & 1.21e-07 \\
        \hline
        \hline
        5 & 256 & 3.82e-05 & 8.94e-05 \\
        5 & 576 & 1.17e-06 & 2.95e-06 \\
        5 & 1600 & 3.27e-08 & 1.02e-07 \\
        5 & 5184 & 1.90e-09 & 6.08e-09 \\
        \hline
        \hline
        6 & 324 & 1.11e-05 & 2.87e-05 \\
        6 & 676 & 1.50e-07 & 4.70e-07 \\
        6 & 1764 & 1.22e-08 & 3.05e-08 \\
        6 & 5476 & 2.98e-10 & 8.90e-10 \\
        \hline
    \end{tabular}
    \caption{Errors for Euler equation for B-Splines ($p$ field).}
    \label{tab:euler_sp}
\end{table}

\clearpage

\section{Additional studies}

The results presented previously have been obtained by using different representations for the solution, but also different kinds of sampling (LHS for networks vs Greville points for B-Splines) and different arithmetic precisions (single for networks vs double for B-Splines). Therefore additional investigations are carried out, before concluding about the role of the representation.

\subsection{Sampling}

\begin{figure}[!ht]
  \begin{subfigure}[b]{0.49\linewidth}
    \centering
    \includegraphics[width=\textwidth]{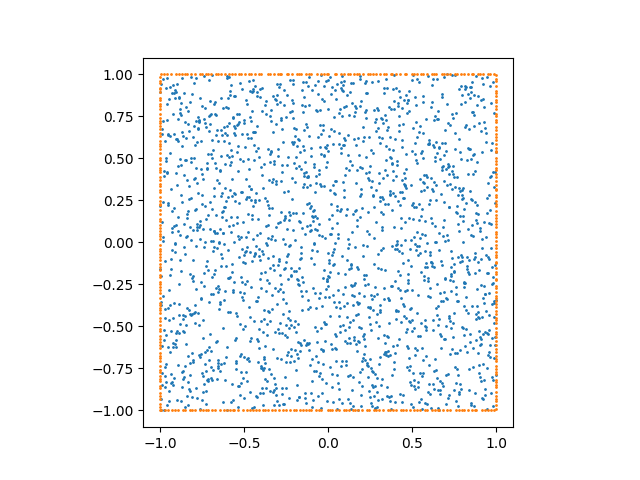}
     \caption{Latin Hypercube Sampling}
    \label{fig:sampling_lhs}
  \end{subfigure}
  \begin{subfigure}[b]{0.49\linewidth}
    \centering
    \includegraphics[width=\textwidth]{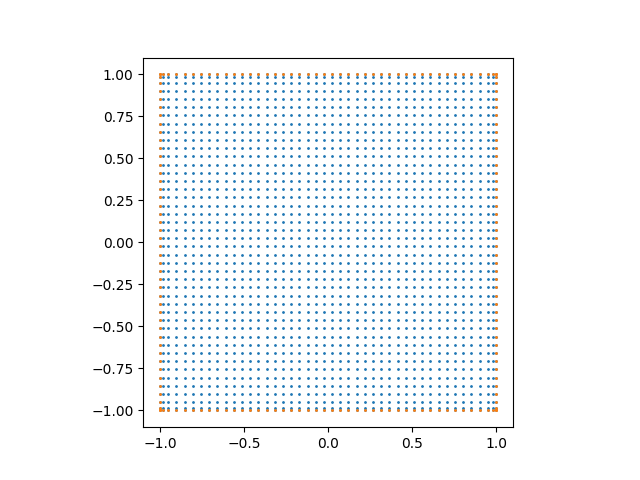}
    \caption{Sampling based on Greville points}
    \label{fig:sampling_greville}
  \end{subfigure}
  \caption{Comparison of sampling methods for an approximate size of 2000 points.}
  \label{fig:sampling}
\end{figure}

Numerical experiments for B-Splines have been achieved by using Greville points to evaluate the residuals. This sampling is derived from a tensor product of distributions and its size increases with the number of degrees of freedom. On the contrary, a random sampling of fixed size has been employed for neural networks. To quantify the impact of this choice on the results, the problem corresponding to the Poisson equation is solved again using neural networks, but employing a sampling similar to Greville points, whose size is chosen to have approximately the same number of networks parameters and collocation points. The difference between the two samplings is illustrated in Fig.~(\ref{fig:sampling}).

\begin{figure}[!ht]
    \centering
    \includegraphics[width=0.5\textwidth]{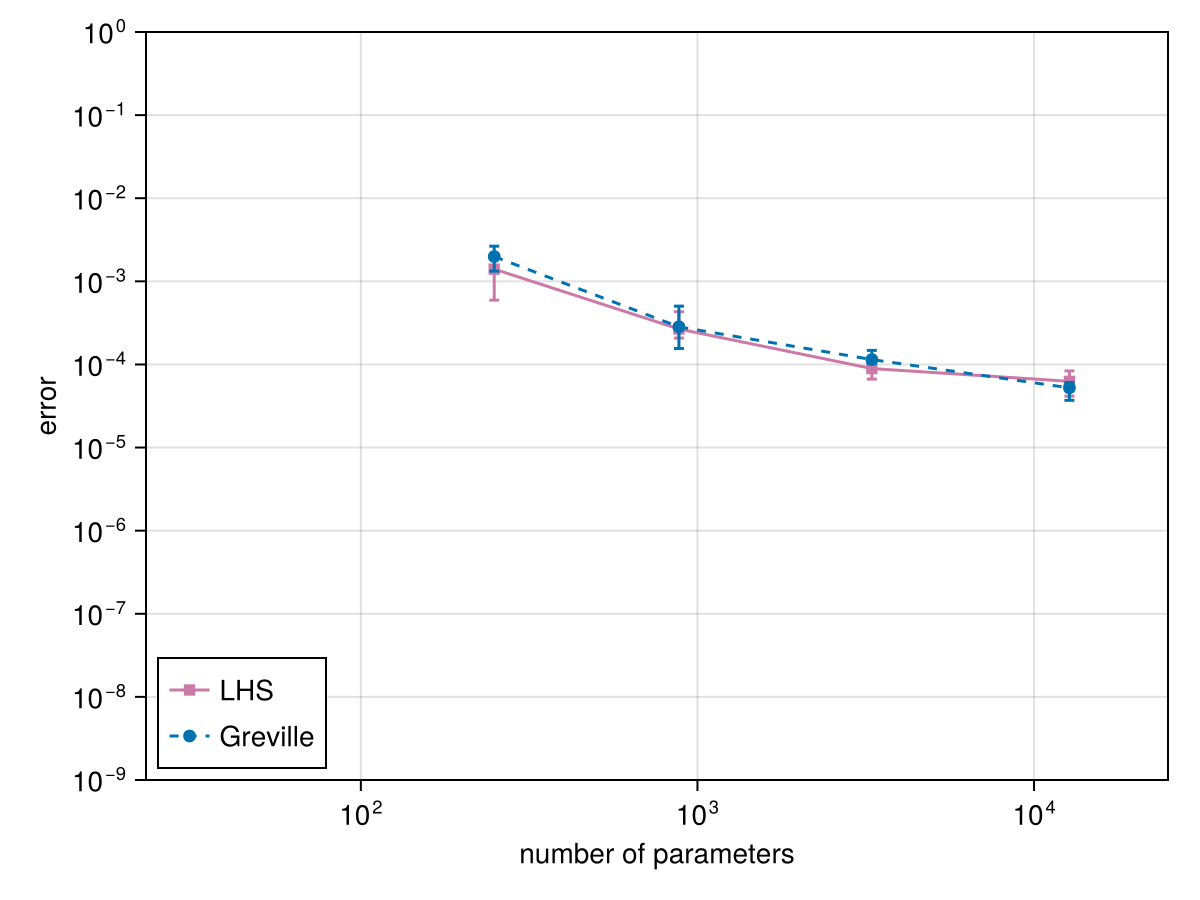}
     \caption{$\epsilon_2$ norm of error w.r.t the number of degrees of freedom for Greville sampling.}
    \label{fig:poisson2d_nn_greville}
\end{figure}

This study is conducted for a set of four-layer networks, including a number of neurons per layer ranging from 8 to 64 as previously. The variation in error as a function of the number of parameters is depicted in Fig.~(\ref{fig:poisson2d_nn_greville}) in the same log scale as Fig.~(\ref{fig:poisson_l2}). Values of error are provided in Tab.~(\ref{tab:sampling}). As seen, the change of sampling type and size has a negligible impact on the results. One can also notice that the variability of the results remains similar, although a deterministic sampling is employed. It proves that this variability mainly originates from the initialization of the networks.

\begin{table}
    \centering
    \begin{tabular}{|c|c|c||c|c|c|c|c|c| }
        \hline
        & Neurons & d.o.f. & \multicolumn{3}{|c|}{$\epsilon_2$} & \multicolumn{3}{|c|}{$\epsilon_\infty$} \\
        \hline
        & & & min & average & max & min & average & max \\
        \hline
        \hline
        LHS & 8 & 249 & 5.94e-04 & 1.40e-03 & 1.98e-03 & 4.09e-03 & 1.19e-02 & 1.83e-02 \\
        & 16 & 881 & 2.08e-04 & 2.67e-04 & 4.31e-04 & 1.93e-03 & 2.88e-03 & 5.06e-03 \\
        & 32 & 3297 & 6.66e-05 & 8.92e-05 & 1.15e-04 & 7.02e-04 & 9.32e-04 & 1.24e-03 \\
        & 64 & 12737 & 4.16e-05 & 6.26e-05 & 8.36e-05 & 5.26e-04 & 7.23e-04 & 1.09e-03 \\
        \hline
        \hline
        Greville & 8 & 249 & 1.33e-03 & 1.99e-03 & 2.65e-03 & 1.07e-02 & 1.26e-02 & 1.42e-02 \\
        & 16 & 881 & 1.56e-04 & 2.83e-04 & 5.03e-04 & 1.39e-03 & 2.56e-03 & 4.64e-03 \\
        & 32 & 3297 & 1.02e-04 & 1.14e-04 & 1.48e-04 & 1.02e-03 & 1.37e-03 & 2.07e-03 \\
        & 64 & 12737 & 3.70e-05 & 5.23e-05 & 6.07e-05 & 4.22e-04 & 5.36e-04 & 7.27e-04 \\
        \hline
    \end{tabular}
    \caption{Errors for Poisson equation using different samplings (four layers).}
    \label{tab:sampling}
\end{table}

\clearpage

\subsection{Single / double precision}

\begin{figure}[!ht]
    \centering
    \includegraphics[width=0.5\textwidth]{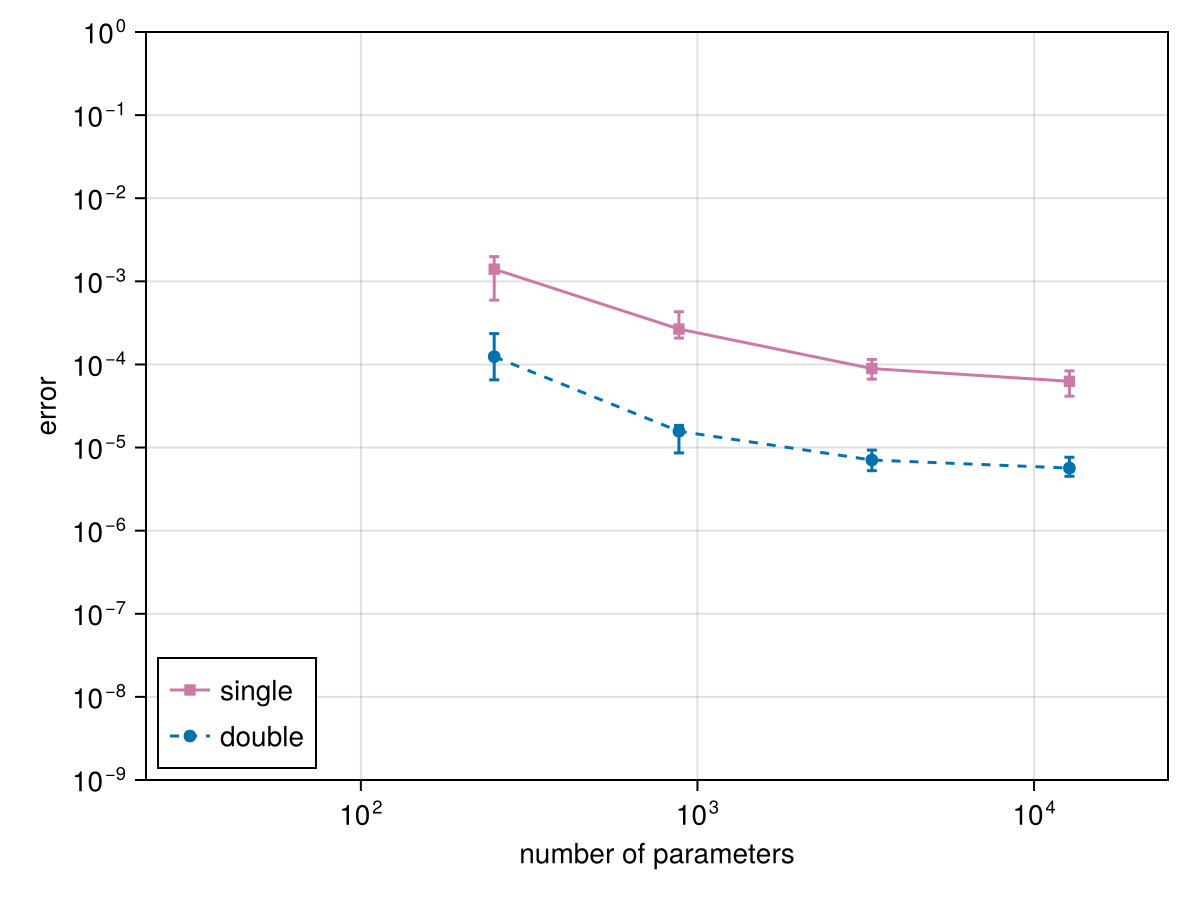}
     \caption{$\epsilon_2$ norm of error w.r.t the number of degrees of freedom for single and double precision.}
    \label{fig:poisson2d_nn_double}
\end{figure}

We investigate now the influence of the arithmetic precision on the results. Indeed, neural networks are usually implemented using single precision, for computational efficiency reeasons, but this choice may have a significant impact, as reported in~\cite{xu2025fp64needrethinkingfailure}. In particular, the use of BFGS-like optimizers during the training may require a higher arithmetic accuracy when upgrading the Hessian matrix from differences of loss gradients. Therefore, the case concerning the Poisson equation is solved again, for the 4-layer networks, using double precision. The variation in error as a function of the number of network parameters is depicted in Fig.~(\ref{fig:poisson2d_nn_double}) and values are provided in Tab.~(\ref{tab:double}). As can be observed, this parameter has a real impact on the results, since the error level drops of one order of magnitude. Nevertheless, the rate of error reduction remains limited: the increase from 32 to 64 neurons per layer yields a small error reduction although a large increase of the number of parameters. Moreover, the computational time has increased significantly.

\begin{table}[!ht]
    \centering
    \begin{tabular}{|c|c|c||c|c|c|c|c|c| }
        \hline
        & Neurons & d.o.f. & \multicolumn{3}{|c|}{$\epsilon_2$} & \multicolumn{3}{|c|}{$\epsilon_\infty$} \\
        \hline
        & & & min & average & max & min & average & max \\
        \hline
        \hline
        single & 8 & 249 & 5.94e-04 & 1.40e-03 & 1.98e-03 & 4.09e-03 & 1.19e-02 & 1.83e-02 \\
        & 16 & 881 & 2.08e-04 & 2.67e-04 & 4.31e-04 & 1.93e-03 & 2.88e-03 & 5.06e-03 \\
        & 32 & 3297 & 6.66e-05 & 8.92e-05 & 1.15e-04 & 7.02e-04 & 9.32e-04 & 1.24e-03 \\
        & 64 & 12737 & 4.16e-05 & 6.26e-05 & 8.36e-05 & 5.26e-04 & 7.23e-04 & 1.09e-03 \\
        \hline
        \hline
        double & 8 & 249 & 6.53e-05 & 1.24e-04 & 2.35e-04 & 7.17e-04 & 1.42e-03 & 2.47e-03 \\
        & 16 & 881 & 8.62e-06 & 1.57e-05 & 1.85e-05 & 9.61e-05 & 2.10e-04 & 3.16e-04 \\
        & 32 & 3297 & 5.27e-06 & 7.08e-06 & 9.32e-06 & 6.51e-05 & 8.41e-05 & 1.14e-04 \\
        & 64 & 12737 & 4.51e-06 & 5.66e-06 & 7.65e-06 & 4.11e-05 & 7.13e-05 & 8.94e-05 \\
        \hline
    \end{tabular}
    \caption{Errors for Poisson equation using single and double precision (four layers).}
    \label{tab:double}
\end{table}

\clearpage

\subsection{Comparison of training phase}

\begin{figure}[!ht]
  \begin{subfigure}[b]{0.49\linewidth}
    \centering
    \includegraphics[width=\textwidth]{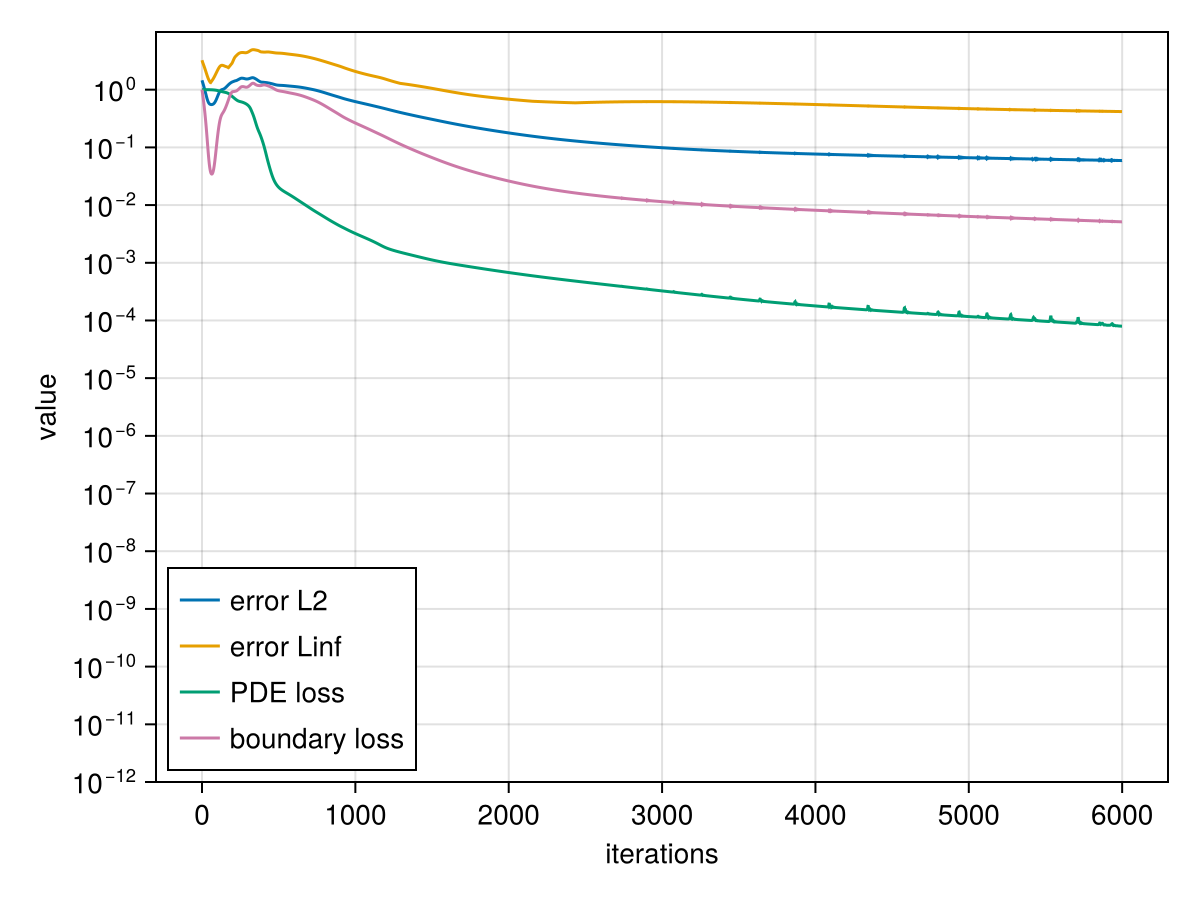}
     \caption{Network using Adam}
    \label{fig:poisson2d_nn_adam}
  \end{subfigure}
  \begin{subfigure}[b]{0.49\linewidth}
    \centering
    \includegraphics[width=\textwidth]{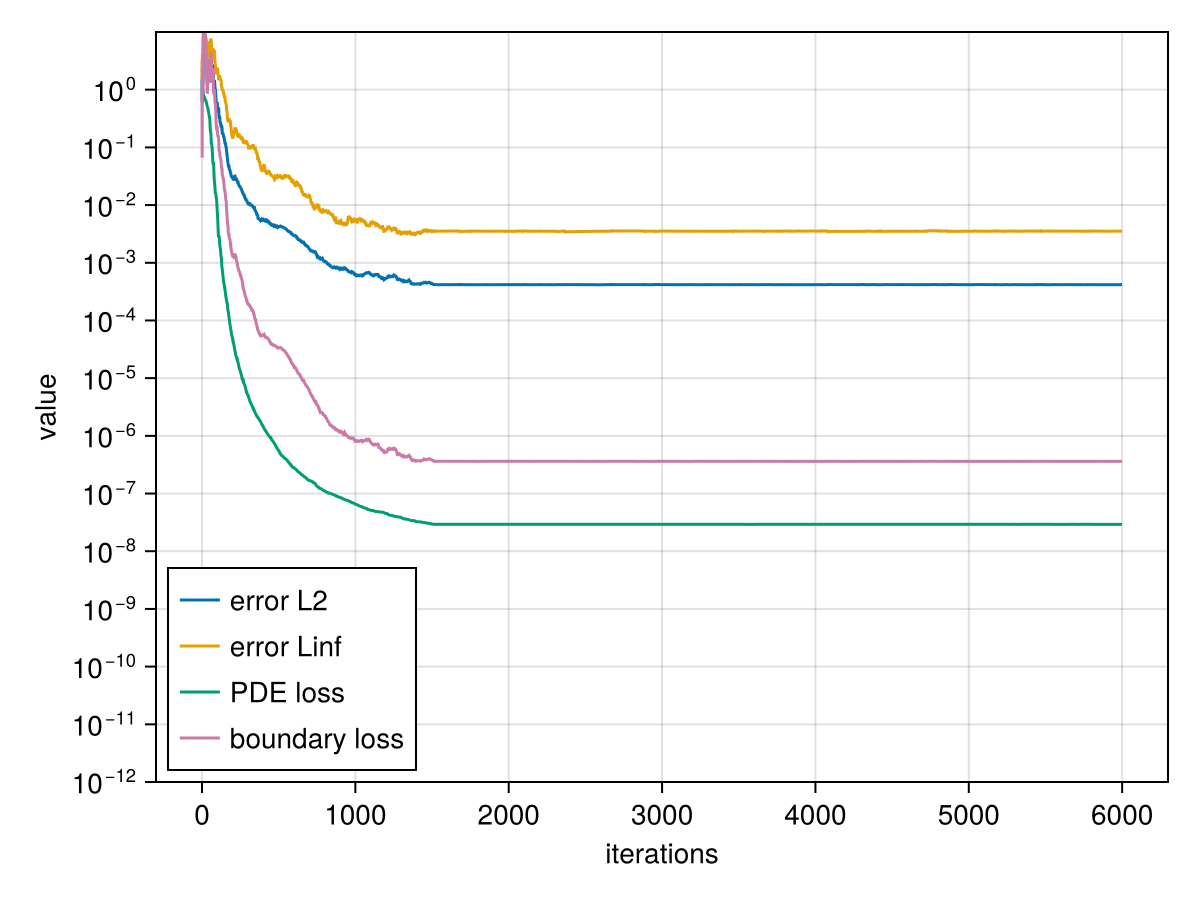}
    \caption{Network using BFGS}
    \label{fig:poisson2d_nn_bfgs}
  \end{subfigure}\\

  \begin{subfigure}[b]{0.49\linewidth}
    \centering
    \includegraphics[width=\textwidth]{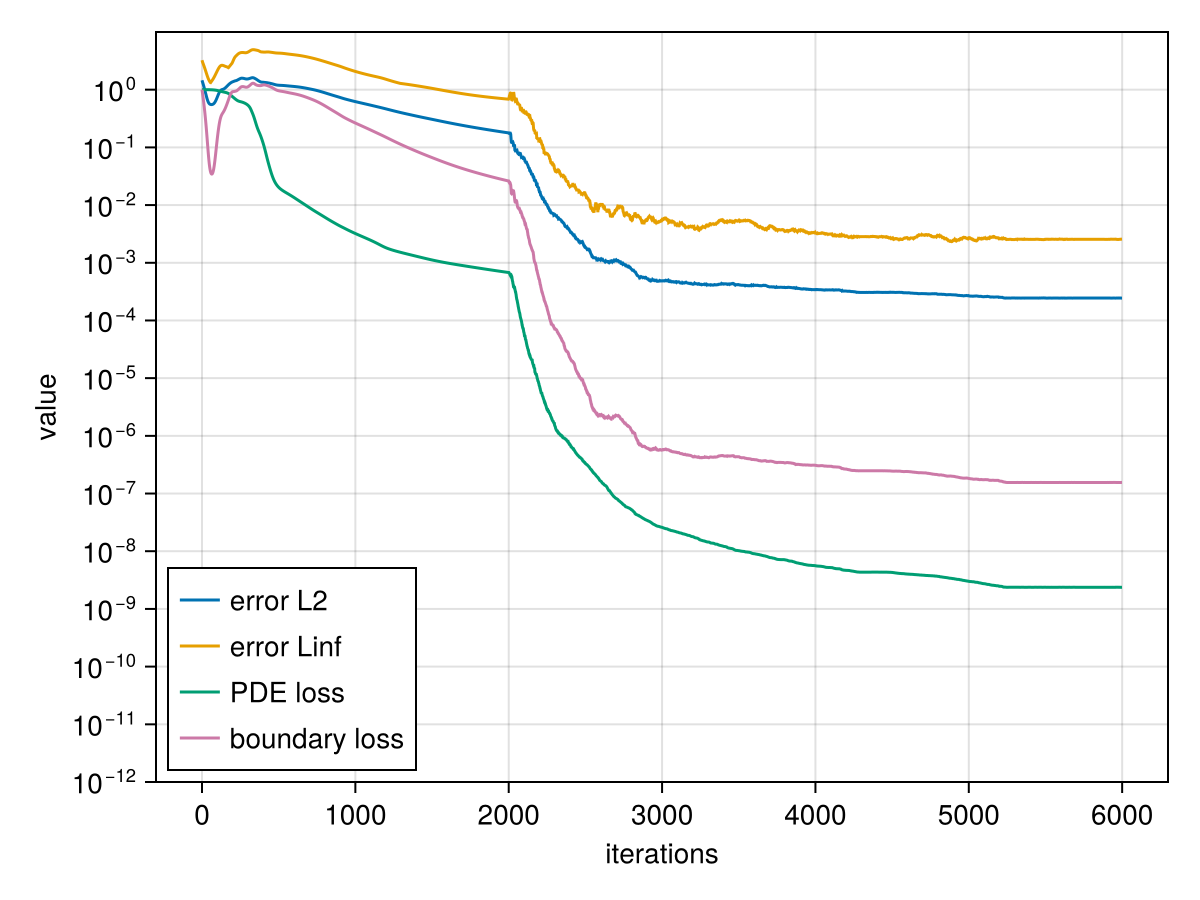}
     \caption{Network using Adam + BFGS}
    \label{fig:poisson2d_nn_adam+bfgs}
  \end{subfigure}
  \begin{subfigure}[b]{0.49\linewidth}
    \centering
    \includegraphics[width=\textwidth]{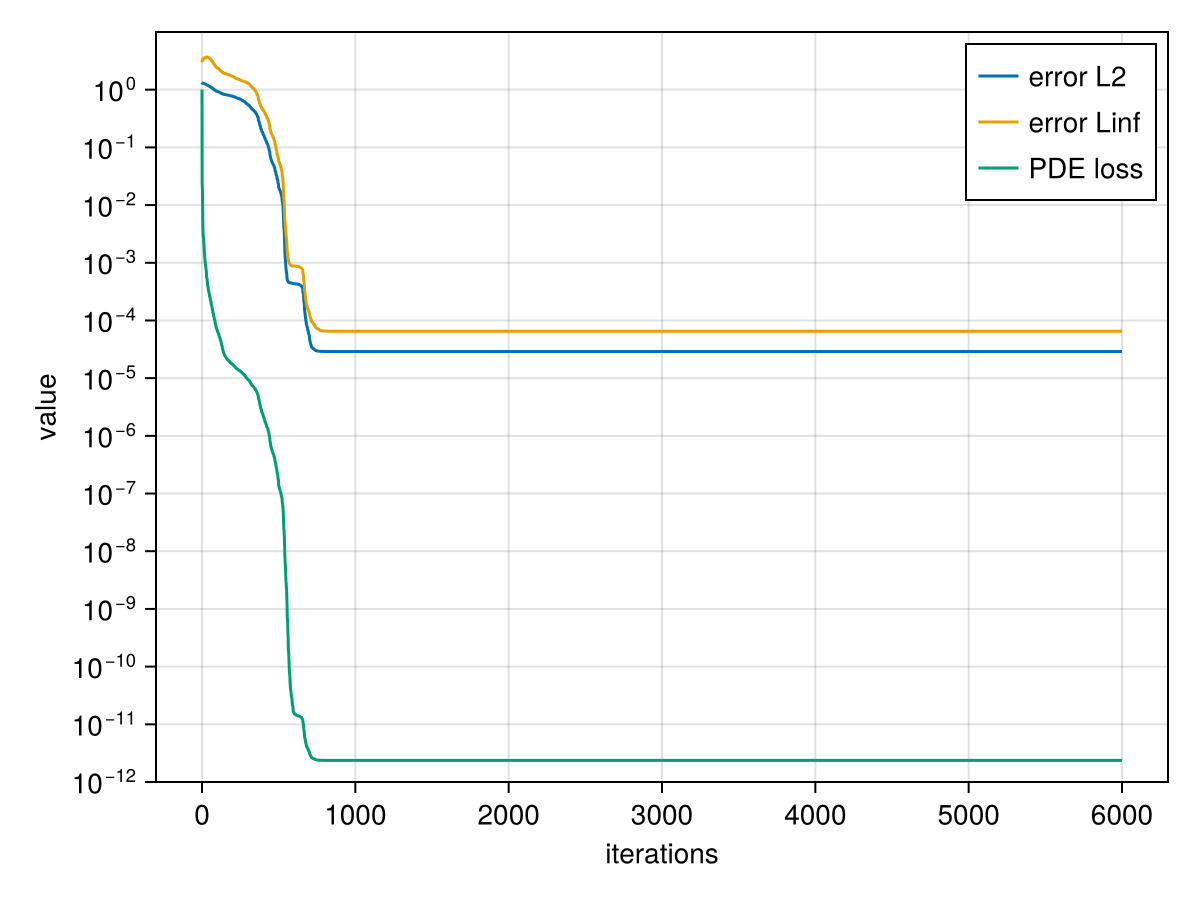}
    \caption{B-Splines using BFGS}
    \label{fig:poisson2d_sp_bfgs}
  \end{subfigure}
  \caption{History of errors and losses during training phase.}
  \label{fig:poisson_optimizer}
\end{figure}

Finally, it is interesting to examine more in depth the minimization of the training loss. As underlined previously, the convexity of the formulation based on B-Splines allows to solve the problem using directly the BFGS algorithm, whereas a first phase using ADAM is required when using neural networks. We illustrate in this section the differences observed using the two optimizers and the two representations. The configuration chosen corresponds to the Poisson case, solved using a 2-layer network with 32 neurons per layer and a B-Spline basis of degree 4 with $34 \times 34$ degrees of freedom, yielding about 1200 parameters for both representations. The history of the training loss functions, $\epsilon_2$ and $\epsilon_\infty$ error norms are plotted in Fig.~(\ref{fig:poisson_optimizer}) using the same scales. Loss values are normalized to be unitary at start.

One can observe a very slow convergence of ADAM optimizer, when used alone (Fig.~(\ref{fig:poisson2d_nn_adam})). If BFGS optimizer is employed directly from the initial network parameters, a fast convergence is obtained (Fig.~(\ref{fig:poisson2d_nn_bfgs})), but towards a local minimum. Indeed, using BFGS after 2000 ADAM iterations yields lower loss values (Fig.~(\ref{fig:poisson2d_nn_adam+bfgs})). However, a far lower loss value is obtained for the B-Spline representation, with fewer iterations. Although these results are obtained for a specific random initialization of the network, very similar plots are generally observed for other initial states, which highlights how the learning problem is conditioned differently depending on the representation used.

\section{Discussion}
\label{sec:disc}

The following lessons can be learnt from these numerical experiments:

\begin{itemize}
    \item Neural networks have the ability to \emph{approximate the solution of a PDE problem using just a few parameters}. This has been noticed in particular when solving the Burgers equation, whose solution is characterized by a localized variation of value. The absence of compact support is certainly an advantage in this context, permitting a self-adaption mechanism focusing the degrees of freedom where needed.
    \item Contrary to B-Splines,\emph{neural networks fail to reduce the error as much as would be desired}. In particular, an increase of the number of parameters does not ensure that the error will be reduced accordingly. No asymptotic convergence has been established, whatever the case.
    \item A significant variability in the results using neural networks has been reported, that seems to be of higher amplitude for convection-dominant problems. This variability is \emph{due to the initialization of the network parameters}, as shown by using a deterministic sampling. The $\epsilon_2$-norm of error typically exhibits a variability of factor 2 to 3, whereas variations of one order of magnitude can be observed for the $\epsilon_\infty$-norm. This is a strong practical limitation of the use of neural networks, since one does not want to perform several random runs when solving an engineering problem.
    \item The resolution of the \emph{loss minimization problem is far more difficult using neural networks} than using B-Splines. It has been shown that the optimization often lead to a local minimum and that the convergence rate may be very low, even if a quasi-Newton optimizer is employed. As consequences, a long iterative process has to be achieved and deceptive results may be obtained. Consequently, it seems that \emph{the error due to the optimization dominates the errors of quadrature and approximation}.
    \item The use of \emph{double precision can help} to increase the efficiency of the training for neural networks, typically reducing the error by one order of magnitude but maintaining a low rate of convergence.
    \item The \emph{use of B-Spline bases yields an efficient physics-informed learning} procedure, including a fast convergence of the loss minimization and an established convergence rate as the number of degrees of freedom increases. This proves that the dominant error for B-Splines is the approximation error. Nevertheless, \emph{the tensorial nature of B-Splines makes this representation untractable} when the dimension of the input (coordinates and other parameters) or the output (physical fields) increases.
\end{itemize}

According to these observations, the bottleneck for physics-informed neural networks lies in the non-convexity and ill-conditioning of the loss function. This yields a slow convergence and a dependancy on the initialization. B-Splines are not affected by these limitations and can exhibit asymptotic convergence rates in accordance with the approximation theory, but suffer from the curse of dimensionality.
This lead to the two following questions: how to improve the convergence of the training of neural networks ? How to improve the scalability of B-Spline bases ?

Regarding neural networks, several works can be found in the literature that attempt to identify the key factors limiting the convergence and find a remedy (see for instance~\cite{krishnapriyan2021characterizing,Pratik_24,Wang_22}). As pointed out in this study, there are certainly several factors increasing the stiffness of the training: the aggregation of different terms in the loss function, the mix of data values and differential information, the non-localization of degrees of freedom, etc.
Several approaches have been proposed to overcome these difficulties, such as adaptive penalty~\cite{Gao_25}, natural gradient descent~\cite{Munk_23}, preconditioning~\cite{Deryck_24}, implicitation of boundary data~\cite{Sukumar_22}, adaptive sampling~\cite{WU2023115671}, etc.
However, as far as we know, none of the proposed methods has permitted to establish a clear convergence as the number of network parameters increases. This property is fundamental for engineering applications because it allows to pave the way towards the control of simulation error.

Concerning B-Splines, an improvement of scalability and computational efficiency could be achieved by introducing either local refinement strategies, such as LR-B-Splines~\cite{BRESSAN2013778} or TH-B-Splines~\cite{GIANNELLI2012485}, that allow to focus the degrees of freedom in regions of interest, and adopt low-rank B-Spline representations~\cite{GEORGIEVA201780,MANTZAFLARIS20171062}, that may yield a reduction of the dimensionality of the problem.

\section{Conclusion}

The objective of this work was to conduct a rigorous comparison between neural and spline representations in the perspective of physics-informed learning. Four PDE problems have been considered and the results have been compared, in terms of $l_2$ and $l_\infty$ error norms, with a focus on the error reduction as the number of degrees of freedom increases. Moreover, the roles of the sampling, the arithmetic precision and the optimizer have been studied more in depth for one case.

It has been found that the two representations exhibit different properties. On the one hand, neural networks can provide a rather good approximation of the solution with a few parameters, but suffers from difficulties during the training that prevent from achieving a high accuracy. On the other hand, B-Splines allow a rigourous control of errors, but have a low scalability due to their tensorial nature. Therefore, a possible perspective would consist in the hybridization or coupling of both representations to benefit from their respective advantages.

\clearpage

\bibliographystyle{plain}
\bibliography{references}

\end{document}